\documentstyle[12pt,bezier,macros,size]{amsart}
\title{Irreducible Constant Mean Curvature 1 Surfaces in\\
        Hyperbolic Space with Positive Genus}
\author{Wayne Rossman, Masaaki Umehara and Kotaro Yamada}
\thanks{
        This research was 
        supported in part by Grant-in-Aid for Scientific Research,
        the Ministry of Education, Science, Sports and Culture, Japan,
        and by a fellowship from the Japan
        Society for the Promotion of Science.}
\begin{document}
\maketitle
\begin{abstract}
 In this work we give a method for constructing a one-parameter family 
 of complete \cmcone\ (i.e.~constant mean  curvature 1)
 surfaces  in hyperbolic 3-space
 that correspond to a given  complete minimal surface with finite total 
 curvature in Euclidean 3-space.   
 We show that this one-parameter family of surfaces with the same symmetry
 properties  exists for all given minimal surfaces satisfying
 certain conditions. The surfaces we construct in this paper 
 are irreducible,
 and in the  process of showing this,  we also prove
 some results about the reducibility of surfaces.
 
 Furthermore, in the case that the surfaces are of genus $0$, we are  able to
 make some estimates on the range of the parameter for the  one-parameter
 family.  
\end{abstract}

\section{Introduction}\label{sec:intro} 

Recently, 
new examples of immersed \cmcone\ surfaces of finite
total curvature in the hyperbolic $3$-space $\hyp$
of curvature $-1$ have been found
(cf.\ \cite{UY1}, \cite{UY2}, \cite{Small}, and \cite{UY4}).
One method used to make these new examples is the following:
The set of all conformal branched \cmcc\ (i.e.\ constant mean 
curvature $c$) immersions in $\hypc$ 
of finite total curvature with hyperbolic Gauss map $G$
defined on a compact Riemann surface $\overline M$ 
corresponds bijectively to the set of 
conformal pseudometrics of constant curvature $1$ with conical 
singularities on $\overline M$. (cf.~\cite{UY4}.)
By the work of Small~\cite{Small}, 
this correspondence can be explicitly written 
when the immersion can be lifted to a null curve in 
$PSL(2,\C)=SL(2,\C)/\{\pm 1\}$. 
This gives a method for constructing new examples.  
However, to construct non-branched 
\cmcc\ surfaces is still difficult, 
because the method above does not give any control 
over branch points.  

In this paper we use a new method to construct new examples 
without branch points, which have higher genus,
many symmetries, and embedded ends. 
More precisely, we prove that for each complete symmetric 
finite-total-curvature minimal surface in $\R^3$ with 
a non-degenerate period problem, 
there exists a corresponding one-parameter family of
\cmcone\  surfaces in $\hyp$.  
We define the terms ``symmetric'' and ``non-degenerate'' later.  
To prove the existence of these corresponding one-parameter families, 
we begin by using a small deformation from the original 
minimal surface in $\R^3$,
preserving its (hyperbolic) Gauss map
$G$ and Hopf differential $Q$.  
This gives us \cmcc\ surfaces in $\hypc$, for $c \approx 0$.  
Finally, we rescale the surfaces into \cmcone\  surfaces in $\hyp$.  
The method is somewhat similar to that of \cite{UY2},
the main difference being that
we use the duality on \cmcone\  surfaces to keep the
symmetry properties of the initial minimal surfaces.

Here, we briefly outline the construction: 
Let  $f_0\colon{}M\to \R^3$ be a conformal minimal immersion 
defined on a Riemann surface $M$. 
We set
\[
    G  =  \frac{\partial x_3}{\partial x_1-i\partial x_2} \, ,
         \quad
    \omega  =  \partial x_1-i\partial x_2 \,,
         \quad
    Q   = \omega \cdot dG ,
\]
where the dot means the symmetric product.
Then $G$ is the Gauss map of $f_0$, and $Q$ is the 
Hopf differential, namely the $(2,0)$-part
of the complexification of the second fundamental form.
The pair $(G,\omega)$ is called the Weierstrass data of $f_0$.
For the \cmcc\ surfaces in $\hypc$, the Weierstrass data can 
also be defined.
The \cmcc\ surface $f_c$ with the same Weierstrass data 
$(G,-\omega)$ as $-f_0$ is only defined on 
the universal cover of $\widetilde M$.
We set 
\[
     \D_{M}^{(c)}(G,Q)
         =\left\{f\colon{}\widetilde M\to\hypc\,;\,
            \begin{array}{l}
              \mbox{$f$ is a conformal \cmcc\ immersion}\\
              \mbox{with hyperbolic Gauss map $G$}\\
              \mbox{and Hopf differential $Q$.}
            \end{array}
          \right\}  .
\]
In Section~\ref{sec:red}, we show that the set
$\D_{M}^{(c)}(G,Q)$ can be identified with the hyperbolic $3$-space 
$\H^3$.  
  (We will use two different notations for hyperbolic space.  
   The notation $H^3$ will be used to represent the ambient space 
    for CMC surfaces, 
    and the notation $\H^3$ will be used to represent the parameter 
    space of $\D(G,Q)$.)
We show that any two \cmcc\ immersions in 
$\D^{(c)}_{M}(G, Q)$ are non-congruent, 
and 
they are dual surfaces of the \cmcc\ surface $f_c$.
  (As mentioned above, $(G,-\omega)$ is the Weierstrass data of $f_c$.
   The duality exchanges the roles of the hyperbolic Gauss map and
   the secondary Gauss map, and at the same time 
   the sign of the Hopf differential is reversed.)

Then, defining a subset $I_{M}^{(c)}(G,Q)$
of immersions which are single-valued on $M$:
\begin{equation}\label{eq:defI}
   I_{M}^{(c)}(G,Q)=
   \left\{
     f\in\D_{M}^{(c)}(G,Q)\,;\,
     \text{$f$ is single-valued on $M$ itself.}
   \right\} \, ,
\end{equation}
we show that $I_{M}^{(c)}(G, Q)$ is either empty, or is a 
connected totally geodesic subspace of dimension $0$, 
$1$ or $3$ in  $\H^3$. 
(The two-dimensional case does not occur.)
When $\dim  I_{M}^{(c)}(G,Q)=0$, the unique \cmcc\ immersion
in $I_{M}^{(c)}(G,Q)$ is irreducible in the sense of \cite{UY1}.
For initial minimal surfaces $f_0$, which are symmetric,
and have non-degenerate period problems,
we can construct a one-parameter family 
of \cmcone\ surfaces $f_c:M\to \hypc$ with the same symmetry
properties as $f_0$.
Moreover, if the initial minimal surface
is generated from a fundamental domain
by reflections with respect to three non-parallel planes,
we can show that $I_{M}^{(c)}(G,Q)$ consists of only one point, 
whenever $c$ is sufficiently small. 
This unique surface coincides with the above $f_c$.
If we consider $\hypc$ as the Poincar\'e model of radius $1/|c|$,
then it converges to the initial surface $f_0$ as $c\to 0$.
(See Remark~\ref{rem:converge}.)
The method we have just described can be applied 
to virtually all of the well-known  symmetric minimal surfaces in $\R^3$.
(See Section~\ref{sec:ftc}.) 
When the hyperbolic Gauss map $G$ and the Hopf differential $Q$ 
have certain properties,
our method is valid even without initial minimal surfaces. 
In Section~\ref{sec:fence}, we demonstrate this 
by constructing hyperbolic correspondences of
the catenoid fence and Jorge-Meeks fences (cf.~\cite{Karcher},
\cite{Rossman}).  
Taking limits of these surfaces as $c$ approaches $0$, 
we obtain alternative proofs of the existence of
the corresponding minimal surfaces in $\R^3$.
 
Other useful applications of duality for \cmcone\ surfaces 
can be found in \cite{UY5}.

\section{Preliminaries}\label{sec:pre}
In this paper,  we use the following identification:
\begin{equation}\label{eq:hyp}
   \hypc=\left\{\frac{1}{c}a a^{*}\,;\, a\in SL(2,\C)\right\},
\end{equation}
where $a^{*}={}^t\bar{a}$.  
The complex Lie group $PSL(2,\C)=SL(2,\C)/\{\pm 1\}$ acts isometrically on
$\hypc$ by $a \circ p=a p a^*$, where $a\in PSL(2,\C)$ 
and $p \in \hypc$. 
(See \cite{Bryant}, \cite{UY1}, or \cite{UY2}.)

Let $M$ be a Riemann surface and 
$f\colon{}M\to \hypc$  a complete conformal \cmcc\ immersion 
of finite total curvature.  
Then there is a null holomorphic immersion 
$F\colon{}\widetilde M \to SL(2,\C)$
defined on the universal cover $\widetilde M$ of $M$
such that $f = (1/c) F F^*$.
Such $F$ is uniquely determined up to the ambiguity 
$F b$, for $b\in SU(2)$. 
Define a meromorphic function $G$ by
\begin{equation}\label{eq:hgauss}
    G=\frac {dF_{11}}{dF_{21}}=\frac {dF_{12}}{dF_{22}} \; , 
\end{equation}
where $F=(F_{ij})_{i,j=1,2}$.
Then $G$ is a meromorphic single-valued function on $M$.  
The function $G$ is called the {\it hyperbolic Gauss map\/} of $f$
(cf.\ \cite{Bryant}).
We can set
\begin{equation}\label{eq:wode}
    F^{-1}dF = c 
               \begin{pmatrix}  
                   g & -g^2 \\
                   1 & -g\hphantom{^2}
               \end{pmatrix} \omega \; , 
\end{equation}
where $g$ is a meromorphic function on $\widetilde M$ and 
$\omega$ is a holomorphic $1$-form defined on $\widetilde M$.
We call the pair $(g,\omega)$ the {\it Weierstrass data\/} 
of the \cmcc\ immersion $f$, 
and $g$ the {\it secondary Gauss map\/} of $f$. 
The pair $(g,\omega)$ also have 
an $SU(2)$-ambiguity with respect to that of $F$
(cf.\ Remark~\ref{rem:cong}).
By using the Weierstrass data $(g,\omega)$, 
the first fundamental form $ds^2$ and 
the second fundamental form $\secondff$ are written as
\begin{equation} \label{eq:first}
   ds^2 = (1+|g|^2)^2\,\omega\cdot \bar \omega, \qquad
            \secondff = -Q-\overline Q + c\,ds^2 \; ,
\end{equation}
where 
\begin{equation}\label{eq:hopf}
    Q=\omega\cdot dg.
\end{equation}
Conversely, if the Weierstrass data $(g,\omega)$ on $M$
satisfying the compatibility condition
\begin{equation}\label{eq:compati}
    ds^2 = (1+|g|^2)^2\,\omega\cdot \bar \omega>0  
\end{equation}
is given abstractly,
then there is a \cmcc\ immersion $f\colon{}\widetilde M\to \hypc$
with Weierstrass data $(g, \omega)$.

\begin{remark}\label{rem:cong}
  Let $\tilde f\colon{} \widetilde M\to \hypc $ be another conformal 
  \cmcc\ immersion with Weierstrass data 
  $(\tilde g, \tilde \omega )$. 
  Then $\tilde f$ and $f$ are congruent if and only if 
  their fundamental forms mutually coincide, 
  which, by (\ref{eq:first}), is equivalent to the following condition: 
  \begin{equation}\label{eq:transf}
     \tilde g=\frac{pg-\overline{q}}{qg+\overline{p}} \; , \qquad 
     \tilde \omega=({qg+\overline{p}})^2 \omega \; , 
  \end{equation}
  for some matrix
  \[ 
      b=\begin{pmatrix}
            p & -\bar{q} \\
            q & \hphantom{-} \bar{p} 
        \end{pmatrix}\in SU(2)\; ,
  \]
  where $|p|^2+|q|^2=1$ (cf.\ \cite[(1.6)]{UY1}).
  We call this $SU(2)$-{\it equivalence}.  
  Transforming $(g,\omega)$ into 
  $(\tilde g, \tilde \omega )$ in \eqref{eq:wode} transforms $F$ into 
  $F b^{-1}$.  
  Thus 
  \[
      \tilde{f} = \frac1c(Fb^{-1})(Fb^{-1})^*
                = \frac1cFF^*=f
  \]
  holds.
\end{remark}

The holomorphic quadratic differential $Q$ defined by \eqref{eq:hopf}
is called the {\it Hopf differential\/} of $f$.  
By \eqref{eq:first}, $Q$ is single-valued on $M$.
The hyperbolic Gauss map $G$, the secondary Gauss map $g$ and the
Hopf differential $Q$ satisfy the identity (cf.\ \cite{UY1}, \cite{UY4}) 
\begin{equation}\label{eq:S}
   S(g)-S(G)=2cQ \; , 
\end{equation}
where $S(g)=S_z(g)dz^2$ is the Schwarzian derivative of $g$, 
namely $S_z(g)$ is 
\begin{equation}\label{eq:schwarz}
  S_z(g)=\left(\frac{g''}{g'}\right)'-\frac12\left(\frac{g''}{g'}\right)^2
  \qquad ({}'=d/dz) \; . 
\end{equation}

\begin{definition}\label{def:lift}
  Let $(g,\omega)$ be a Weierstrass data on $M$ with the
  compatibility condition 
  \eqref{eq:compati}. 
  Then there exists a conformal \cmcc\ immersion 
   $f\colon{}\widetilde M\to \hypc$ with the Weierstrass data
   $(g,\omega)$. Let $G$ be the hyperbolic Gauss map of $f$.
  Then there exists a unique null holomorphic immersion
  $F\colon{}\widetilde M\to SL(2,\C)$ satisfying \eqref{eq:hgauss} 
  and \eqref{eq:wode} (cf.~\cite[Theorem 1.6]{UY4}). 
  We call $F$ the {\it lift\/} of $f$ with respect to the 
  Weierstrass data $(g,\omega)$.
\end{definition}

For a meromorphic function $G$ and  a holomorphic 
$2$-differential $Q$ defined on  $M$, 
the set $\D_{M}^{(c)}(G,Q)$ is defined as in the introduction.

\begin{lemma}\label{lem:Dempty}
   The following assertions are equivalent\rom{:}
   \begin{enumerate}
      \item\label{item1:Dempty}
           $\D_{M}^{(c)}(G,Q)$ is nonempty.
      \item\label{item3:Dempty}
           The symmetric tensor
           \[
                ds_G^2= (1+|G|^2)^2\,\frac{Q}{dG}\cdot
                                \overline{\left(\frac{Q}{dG}\right)}
           \]
           is positive definite on $M$.
        \end{enumerate}
  Moreover,  any immersion $f$ in $\D_{M}^{(c)}(G,Q)$ has 
  a  complete induced metric whenever
  $ds_G^2$ is a complete metric of finite total curvature
  and $Q$ has poles of order at most\/ $2$. 
\end{lemma}
\begin{pf}
  Assume that \ref{item3:Dempty} holds.
  Then there is a conformal \cmcc\ immersion
  $f\colon{}\widetilde M \to \hypc$ with the Weierstrass data $(G,-Q/dG)$.
  Let $g$ be the hyperbolic Gauss map of $f$. 
  Take the lift  $F\colon{} \widetilde{M}\rightarrow SL(2,\C)$  of $f$. 
  We set 
  \[
       f^\#=\frac1cF^{-1}(F^{-1})^*\colon{}\widetilde M\to \hypc,
  \]
  which is called the dual surface in \cite{UY5}.
  By taking the dual, the hyperbolic Gauss map $G$ and
  the secondary Gauss map $g$ are exchanged, and  the sign of $Q$ is reversed.
  So the conformal \cmcc\ immersion $f^\#$ has  
  the hyperbolic Gauss map $G$ and the Hopf differential $Q$.
  (See Proposition~4 of \cite{UY5}, 
   with special attention to the sign of $Q$.)
  Hence $f^{\#}\in\D_{ M}^{(c)}(G,Q)$, and 
  \ref{item3:Dempty} implies \ref{item1:Dempty}.  
  Moreover the completeness of the induced metric of $f^\#$
  follows from \cite[Lemma 5]{UY5}.

  Conversely, suppose that $\hat f\in \D_{M}^{(c)}(G,Q)$ is given. 
  Let $(\hat g,\hat \omega)$ be the Weierstrass data of $\hat f$.
  Taking the lift $\hat F$ of $\hat f$ with respect to 
  the Weierstrass data $(\hat g,\hat \omega)$,
  the conformal \cmcone{} immersion defined by
  $\hat f^\#=(\hat F^{-1})(\hat F^{-1})^*$ has the Weierstrass data
  $(G,-Q/dG)$. 
  Since $\hat F^{-1}$ is an immersion, the first fundamental form
  $ds_G^2$ of $\hat f^\#$ is positive definite.
\end{pf}

\begin{corollary}\label{cor:add}
  Suppose that $f\in \D_{M}^{(c)}(G,Q)$ satisfies
  $f=(1/c)FF^*$, where $F\colon{}M\to SL(2,\C)$ is a null holomorphic
  immersion. 
  Then $F$ satisfies the differential equation
  \begin{equation}\label{eq:dw}
      dF\cdot F^{-1}=c
                     \begin{pmatrix}
                        G & -G^2 \\
                        1 & -G\hphantom{^2}  
                     \end{pmatrix}
                      \frac{Q}{dG}. 
  \end{equation}
\end{corollary}
\begin{pf}
  As seen above, 
  the dual immersion $f=(1/c)(F^{-1})(F^{-1})^*$ has the Weierstrass
  data $(G,-Q/dG)$. 
  Since $(F^{-1})^{-1}d(F^{-1})=-dFF^{-1}$, the assertion
  immediately follows from \eqref{eq:wode}.
\end{pf}
\begin{lemma}\label{lem:Dparam}
   If $\D_{M}^{(c)}(G,Q)$ is not empty, 
   then it is identified with the hyperbolic $3$-space $\H^3$.
\end{lemma}

In the above statement, we used the notation $\H^3$
for the hyperbolic $3$-space to distinguish it from the
hyperbolic $3$-space as the ambient space for CMC surfaces.

\begin{pf}
  Choose a \cmcc\ immersion $f_0\in \D_{M}^{(c)}(G,Q)$, 
  and fix a null holomorphic immersion $F_0\colon{}\widetilde M\to SL(2,\C)$ 
  such that $f_0=(1/c)F_0^{}F_0^*$.
  Consider any immersion $f_1\in\D_{M}^{(c)}(G,Q)$.
  Then there exists a null holomorphic immersion
  $F_1\colon{}\widetilde M\to SL(2,\C)$ satisfying $f_1 = (1/c) F_1^{}
F_1^*$.  
  Let $g_0$ and $g_1$ be the secondary Gauss maps of $f_0$ and $f_1$,
  respectively.  
  By  \eqref{eq:S}, $S_z(g_0) = S_z(g_1)$.
  Thus a well-known property of the Schwarzian derivative yields 
  \[
      g_1 = \frac{a_{11}g_0+a_{12}}{a_{21}g_0+a_{22}} 
      \qquad \mbox{for some} \qquad 
      a = (a_{kj})_{k,j=1,2} \in SL(2,\C) \; .
  \]
  Since the hyperbolic Gauss maps of $f_0$ and $f_1$ coincide, 
  we have $F_1=F_0 a^{-1}$ (cf.\ \cite[(1.6)]{UY4}).
  Thus $f_1=(1/c) F_0^{} (a^{-1}) (a^{-1})^* F_0^{*}$, and 
  \[
     \begin{aligned}
       \D_{M}^{(c)}(G,Q)&=
               \left\{
                    \frac{1}{c} F_0 (a^{-1}) (a^{-1})^* F_0^*\,;\,
                    a\in SL(2,\C)
               \right\}\\
          &\simeq
               \left\{ (a^{-1}) (a^{-1})^*\,;\,a\in SL(2,\C)\right\}
                =\H^3.\qquad \Box
    \end{aligned}
  \]
  \renewcommand{\qed}{}
\end{pf}

By Remark~\ref{rem:cong}, the following assertion is immediately obtained.

\begin{corollary}\label{cor:noncong}
  For any two distinct \cmcc\ immersions $f_1,f_2\colon{} \widetilde M\to
\hypc$
  in $\D_{M}^{(c)}(G,$ $Q)$, there is no isometry $T \in
  SL(2,\C)$ such that $T(f_1)=f_2$.
\end{corollary}

\section{Reducibility}\label{sec:red}

Let $M$ be a  Riemann surface and $p\colon{}\widetilde M\to M$ 
the universal cover of $M$.   
We fix a reference point $\tilde z_0\in \widetilde M$ and identify 
canonically the fundamental group $\pi_1(M)=\pi_1(M,p(\tilde z_0))$ 
with the deck transformation group on $\widetilde M$.
Take a meromorphic function $G$ and a holomorphic $2$-differential $Q$ 
defined on $M$.
We identify the lifts $G\circ p$ and $Q\circ p$ on $\widetilde M$ 
with $G$ and $Q$ themselves.

We consider the subset $I_M^{(c)}(G,Q)$ of 
$\D_{M}^{(c)}(G,Q)$ as defined as \eqref{eq:defI}.
When $\D_{M}^{(c)}(G,Q)$ is non-empty,
as seen in Lemma~\ref{lem:Dparam}, 
we can consider $I_{M}^{(c)}(G,Q)$ as a subset of 
the hyperbolic $3$-space $\H^3$.
The set $I_{M}^{(c)}(G,Q)$ is closely related to the reducibility 
of \cmcc\ surfaces. 
We recall the definition (cf.\ \cite[Definition~3.1]{UY1}):
For $f\in \D_{M}^{(c)}(G,Q)$, 
we define a pseudometric $d\sigma^2_f$ by
$d\sigma_f^2=(-K)ds^2$, where $K$ is the sectional
curvature of the first fundamental form $ds^2$. 
Then (cf.\ \cite[(2.8)]{UY4})
\begin{equation}\label{eq:Def}
    d\sigma^2_f = \frac {4\,dg\cdot d\overline{g}}{(1+|g|^2)^2},
\end{equation}
where $g$ is the secondary Gauss map.
By \eqref{eq:first}, \eqref{eq:hopf} and \eqref{eq:Def}, 
one can easily get the relation
\begin{equation}\label{eq:metHop}
    ds^2\cdot d\sigma_f^2=4\,Q\cdot \overline{Q}.
\end{equation}
By \eqref{eq:Def}, 
$d\sigma_f^2$ is the pull back of the canonical metric on
the unit sphere $S^2\cong \C\cup\{\infty\}$ by 
the secondary Gauss map $g$.
For each deck transformation $\tau\in \pi_1(M)$, 
there exists a matrix $\tilde\rho(\tau)=(a_{jk})_{j,k=1,2}$ 
in $PSL(2,\C)$ such that
\begin{equation}\label{eq:pull}
    g\circ \tau^{-1}=\tilde\rho(\tau)*g =
         \frac{a_{11}g+a_{12}}{a_{21}g+a_{22}},
\end{equation}
where the asterisk  means the
action of $\tilde \rho(\tau)$ as a M\"obius transformation. 
Then a representation $\tilde\rho\colon{}\pi_1(M)\to PSL(2,\C)$ is induced.
Let $F\colon{}\widetilde M\to SL(2,\C)$ 
be the lift of $f$ with respect to the Weierstrass data $(g,Q/dg)$ 
(cf. Definition~\ref{def:lift}).  
Since $F\circ \tau$ and $F\tilde \rho(\tau)$ have the same
hyperbolic Gauss map and the secondary Gauss map $g\circ \tau$,
Theorem~1.6 in \cite{UY4} yields
\begin{equation}\label{eq:decktr}
    F\circ \tau =s_\tau F\, \tilde \rho(\tau) \qquad (\tau\in \pi_1(M)),
\end{equation}
where $s_{\tau}=1$ or $-1$.
We set
\[
    \rho(\tau)=s_\tau \tilde \rho(\tau)\qquad 
    (\tau\in \pi_1(M)).
\]
Then $\tilde \rho\colon{}\pi_1(M)\to PSL(2,\C)$ can be lifted to
a representation
\[
    \rho\colon{}\pi_1(M)\to SL(2,\C).
\]
Since the representation $\rho$ depends on the choice of 
the Weierstrass data,
we will also use the notation  $\rho_F^{}$ 
as well as $\rho$  when such  an explicit notation is required.
\begin{lemma}\label{lem:addi}
    Suppose $f\in \D_M^{(c)}(G,Q)$. 
    Then $f\in I_M^{(c)}(G,Q)$ 
    if and only if $\rho(\tau)\in SU(2)$ for all $\tau\in \pi_1(M)$.
\end{lemma}
\begin{pf}
  Since $G$ and $Q$ are single-valued on $M$,
  by Corollary~2.4 in \cite{UY4},
  $f\in I_M^{(c)}(G,Q)$ if and only if $ds^2$ is single-valued on $M$.
  The condition is equivalent to $d\sigma_f^2$ being single-valued on $M$
  by \eqref{eq:metHop}.
  Thus $f\in I_M^{(c)}(G,Q)$ if and only if
  each $\rho(\tau)$ preserves $d\sigma_f^2$, that is,
  $\rho(\tau)\in SU(2)$.
\end{pf}

A \cmcc\ immersion $f\in I_M^{(c)}(G,Q)$ is said to be 
{\it reducible\/} if
$\rho(\tau_1)\rho(\tau_2)=\rho(\tau_2)\rho(\tau_1)$
holds for all $\tau_1,\tau_2 \in \pi_1(M)$.  
In Theorem~3.3 of \cite{UY1}, 
the second and third authors showed that reducible \cmcc\ surfaces 
have a nontrivial deformation associated with the deformation of 
the Weierstrass data   
$(\lambda g, (1/\lambda)\omega)\,\, (\lambda \in \R\setminus\{0\})$.
This deformation preserves the Hopf differential $Q=\omega\cdot dg$.
Moreover, by \eqref{eq:S}, $S(G)$ is also preserved.
Since the Schwarzian derivative $S(G)$ is invariant under  
the M\"obius transformations of $G$, 
we can place each of these one-parameter family of surfaces
by a suitable rigid motion in $\hypc$ so that the hyperbolic Gauss map
is not changed. 
Thus we get a non-trivial deformation in $I_{M}^{(c)}(G,Q)$
for reducible \cmcc\ surfaces.
As a refinement of this observation, 
we prove the following theorem.
\begin{theorem}\label{thm:red}
 Assume that the subset $I_{M}^{(c)}(G,Q)$ 
  of $\D_{M}^{(c)}(G,Q)$ is not empty.
  Then the set $I_{M}^{(c)}(G,Q)$ is a point,
  a geodesic line $\H^1$
or all of $\H^3=\D_{M}^{(c)}(G,Q)$.
  Moreover, each $f\in I_{M}^{(c)}(G,Q)$ is irreducible 
  if and only if $I_{M}^{(c)}(G,Q)$ is a point.
\end{theorem}
\begin{remark}\label{rem:totgeod}
  In other words, $I_{M}^{(c)}(G,Q)$ is a connected totally geodesic
  $n$-dimensional subset of $\H^3$ with $n = 0$, $1$, or $3$.
  In the case $I_{M}^{(c)}(G,Q)=\H^3$,
  any null holomorphic immersion $F\colon{}\widetilde M\to PSL(2,\C)$ 
  with the hyperbolic Gauss map $G$ and the Hopf differential $Q$
  is single-valued on $M$ 
  (cf. \cite[Theorem~1.6]{UY3}).
\end{remark}

We set $\Gamma:=\rho(\pi_1(M))\subset \mbox{SU}(2)$.
Then Theorem~\ref{thm:red} is an immediate consequence of 
the lemma in the appendix.

\section{A representation of the fundamental group}\label{sec:fund}

In this section, we consider how to find a \cmcc\ immersion 
$f\in\D_{M}^{(c)}(G,Q)$ for any given hyperbolic Gauss map $G$ and 
Hopf differential $Q(\not\equiv 0)$ on a Riemann surface generated by
reflections.
As in the previous section, we fix a reference point
$\tilde z_0\in \widetilde M$ and identify canonically the fundamental group 
$\pi_1(M)=\pi_1(M,p(\tilde z_0))$ with the deck transformation group 
on $\widetilde M$.
Let $\tilde\mu_1$, \dots, $\tilde\mu_m$ be reflections of $\widetilde M$,
that is, conformal orientation-reversing involutions on $\widetilde M$.
Suppose that $\tilde \mu_1,...,\tilde \mu_m$ generate the deck
transformation group on $\widetilde M$ and
induce reflections $\mu_1,\dots,\mu_m$ of $M$
such that
\begin{equation}\label{eq:aditional}
       p\circ \tilde \mu_j= \mu_j \circ p 
       \qquad (j=1,\dots,m),
\end{equation}
where $p\colon{}\widetilde M\to M$ is the covering projection.

Take a pair $(G,Q)$ of a meromorphic function and a holomorphic 
$2$-differential on $M$ such that 
$ds_G^2= (1+|G|^2)^2\,(Q/dG)\cdot\overline{\left(Q/dG \right)}$
is positive definite.
In addition, we assume that
\begin{equation}\label{eq:refassum}
    \overline {Q \circ \mu_j}=Q
        \qquad
        \text{and}
        \qquad
    ds^2_G \circ \mu_j = ds^2_G
\end{equation}
hold. 
By Lemma~\ref{lem:Dempty}, the set $\D_{M}^{(c)}(G,Q)$ is non-empty.

\begin{lemma}\label{lem:gaussref}
  There exist matrices $\sigma(\mu_j)\in SU(2)$ 
  $(j=1,\dots,m)$ such that
  \[
     \overline{G\circ\mu_j}=\sigma(\mu_j)^{-1}* G = 
           \frac{a_{11}^{(j)}G+a_{12}^{(j)}}
     {a_{21}^{(j)}G+a_{22}^{(j)}} \; ,
  \]
  where $\sigma(\mu_j)^{-1}=(a_{k,l}^{(j)})$.
\end{lemma}
\begin{pf}
  Consider a pseudometric
  $d\rho_G^2=4\,dG\cdot d\overline G/(1+|G|^2)^2$
  on $\widetilde M$.
  By \eqref{eq:refassum} and
  \eqref{eq:metHop}, $d\rho_G^2\circ\mu_j=d\rho_G^2$ holds.
  This occurs if and only if there exists a matrix 
  $\sigma(\mu_j)\in SU(2)$ such that 
  $\overline{G\circ\mu_j}=\sigma(\mu_j)^{-1}* G$,
  since $\mu_j$ is orientation reversing.
\end{pf}                
    
Each matrix $\sigma(\mu_j)$~$(j=1,\dots,m)$ is determined
up to sign. 
From now on, we fix a sign of $\sigma(\mu_j)$.

\begin{lemma}\label{lem:repref}
  Let $f\in\D_{M}^{(c)}(G,Q)$ and $F\colon{}\widetilde M\to SL(2,\C)$ 
  be a null holomorphic immersion such that
  $f=(1/c)FF^*$.
  Then there exists a unique matrix 
  $\hat\rho^{}_F(\tilde\mu_j)\in SL(2,\C)$ 
  \rom($j=1,\dots,m$\rom) such that
  \begin{equation}\label{eq:lepref}
    \overline{F\circ\tilde\mu_j}
        =\sigma(\mu_j)^{-1} F \hat\rho^{}_F(\tilde\mu_j)
  \end{equation}
  holds for $\sigma(\mu_j)$ in Lemma~\ref{lem:gaussref}.
\end{lemma}
\begin{pf}
  Consider two \cmcone{} immersions defined by
  \[
       f_1:=(F^{-1})(F^{-1})^*, \qquad
       f_2:=\left({\overline{F\circ\tilde\mu_j}}
                   \vphantom{ 
                     \overline{(F\circ\tilde\mu_j)}^{-1}
                   }
            \right)^{-1}
            \left(\overline{(F\circ\tilde\mu_j)}^{-1}\right)^*.
  \]
  Then by assumption \eqref{eq:refassum}, 
  the two immersions have the same first fundamental form $ds_G^2$ and
  the same second fundamental form $-Q-\overline Q+c\,ds_G^2$.
  Hence these immersions are congruent by the fundamental theorem for
surfaces,
  and so there is a unique isometry $a\in SL(2,\C)$ such that
  $f_1=a f_2 a^*$.
  In particular, there exists $b\in SU(2)$ such that
  (cf. Remark~\ref{rem:cong})
  $F^{-1}=a(\overline{F\circ\tilde\mu_j})^{-1}b$, 
  that is,
  $\overline{F\circ\tilde\mu_j}=b F a$.
  Since $\overline{F\circ\tilde\mu_j}$ has the hyperbolic Gauss map 
  $\overline{G\circ \mu_j}$,
  we have $b=\varepsilon\sigma(\mu_j)^{-1}$ ($\varepsilon=\pm 1$) 
  by Lemma~\ref{lem:gaussref}
  and by \cite[Theorem~1.6]{UY4}.
  Finally, if we set $\hat\rho_F(\tilde \mu_j)=\varepsilon a$, we have the
  desired expression.
  The uniqueness of  $\hat\rho_F(\tilde \mu_j)$ 
  follows from that of the matrix $a$ and the sign $\varepsilon$.
\end{pf}

\begin{remark}\label{rem:refrel}
  For any two reflections $\tilde\mu_j$, $\tilde\mu_k$, 
  the following relation holds:
  \[
      F\circ\tilde\mu_j\circ\tilde\mu_k=
          \overline{\sigma(\mu_k)^{-1}(F\circ\tilde\mu_j)
                    \hat\rho^{}_{F}(\tilde\mu_k)}
          \left(\sigma(\mu_j)
                \overline{\sigma(\mu_k)}
          \right)^{-1} F
          \hat\rho^{}_{F}(\tilde\mu_j) 
          \overline{\hat\rho^{}_{F}(\tilde\mu_k)}.
  \]
\end{remark}

\begin{corollary}\label{lem:repchange}
  Let $f\in\D_{M}^{(c)}(G,Q)$ and $F\colon{}\widetilde M\to SL(2,\C)$ 
  be a null holomorphic immersion such that
  $f=(1/c)FF^*$. 
  Then 
  $\hat\rho^{}_{Fa}(\mu_j)=a^{-1} \hat\rho^{}_F(\mu_j)\bar a$
  holds for any $a\in SL(2,\C)$.
\end{corollary}

\begin{lemma}\label{lem:id}
  $\sigma(\mu_j)\overline{\sigma(\mu_j)}=\id$.
\end{lemma}

\begin{pf}
  Since $\tilde\mu_j\circ\tilde\mu_j=\id$, 
  $G=G\circ\mu_j\circ \mu_j=
   \{\sigma(\mu_j)  \overline{\sigma(\mu_j)}\}^{-1}* G$,
  and we see that $\sigma(\mu_j)\overline{\sigma(\mu_j)}=\pm\id$.
  Suppose that $\sigma(\mu_j)\overline{\sigma(\mu_j)}=-\id$.
  Since $\sigma(\mu_j)\in SU(2)$, we have 
  \[ 
      \sigma(\mu_j) = \pm \begin{pmatrix}
                              \hphantom{-}0   & 1 \\
                              -1   & 0
                          \end{pmatrix} \; \; ,  
  \]
  and so $G\circ \mu_j=\sigma(\mu_j)^{-1}*G=-1/\overline G$.
  Let $z$ be a fixed point of $\mu_j$. Then $|G(z)|^2=-1$ holds,
  a contradiction.
\end{pf}

\begin{remark}\label{rem:ref2}
  By Remark~\ref{rem:refrel} and the previous lemma,
  we have 
  $\hat\rho^{}_F(\mu_j)\overline{\hat\rho^{}_F(\mu_j)}=\id$.  
\end{remark}

Recall that for each $\tau\in \pi_1(M)$, there is a representation
$\rho_F\colon{}\pi_1(M)\to SL(2,\C)$ such that $F\circ \tau=F\rho_F(\tau)$.
(See Section~\ref{sec:red}.)
Each $\tau$ is represented as
\begin{equation}\label{eq:add1}
  \tau=\tilde\mu_{j_1}\circ
       \tilde\mu_{j_2}\circ\dots\circ\tilde\mu_{j_{2k-1}}
       \circ\tilde\mu_{j_{2k}}.
\end{equation}
Then by Remark~\ref{rem:refrel}, we have 
\begin{multline}\label{eq:add2}
  F\circ\tau=
      \left\{{\sigma(\mu_{j_1})}\overline{\sigma(\mu_{j_2})}
              \cdots
              {\sigma(\mu_{j_{2k-1}})}
                \overline{\sigma(\mu_{j_{2k}})}\right\}^{-1}
             \cdot F \cdot \\ 
             {\hat\rho^{}_F(\tilde\mu_{j_1})}
              \overline{\hat\rho^{}_F(\tilde\mu_{j_2})}
              \cdots
             {\hat\rho^{}_F(\tilde\mu_{j_{2k-1}})}
              \overline{\hat\rho^{}_F(\tilde\mu_{j_{2k}})}
                 \; .
\end{multline}
On the other hand, by Lemma~\ref{lem:gaussref},
\[
     G\circ\tau=
     \left\{
      {\sigma(\mu_{j_1})}\overline{\sigma(\mu_{j_2})}
                \cdots
                {\sigma(\mu_{j_{2k-1}})}
              \overline{\sigma(\mu_{j_{2k}})}
       \right\}^{-1}
       * G
\]
holds.
Since $G$ is single-valued on $M$,  we have
\begin{equation}\label{eq:Add}
   {\sigma(\mu_{j_1})}\overline{\sigma(\mu_{j_2})}
   \cdots
   {\sigma(\mu_{j_{2k-1}})}
   \overline{\sigma(\mu_{j_{2k}})}
         =\varepsilon_\tau\id \; ,
\end{equation}        
where $\varepsilon_\tau=1$ or $-1$. 
So we have an expression
\begin{equation}\label{eq:add3}
  \rho^{}_F(\tau)=\varepsilon_\tau
                  \hat\rho^{}_F(\tilde\mu_{j_1})
                  \overline{\hat\rho^{}_F(\tilde\mu_{j_2})}
                  \cdots
                  \hat\rho^{}_F(\tilde\mu_{j_{2k-1}})
                  \overline{ \hat\rho^{}_F (\tilde\mu_{j_{2k}})}.
\end{equation}
The following criterion is useful to construct examples.

\begin{proposition}\label{prop:su2cond}
  Let $f\in\D_{M}^{(c)}(G,Q)$ and $F\colon{}\widetilde M\to SL(2,\C)$ 
  be a null holomorphic immersion such that
  $f=(1/c)F F^*$.
  Suppose
  $\hat\rho^{}_{F}(\tilde\mu_j)\in SU(2)$ \rom($j=1,\dots,m$\rom)
  holds.
  Then the immersion $f$ is an element of $I_{M}^{(c)}(G,Q)$.
  Moreover, each reflection $\mu_j$ extends to an isometry of 
  $\hypc$ preserving the image of $f$.
\end{proposition}
\begin{pf}
  The first assertion immediately follows from Lemma~\ref{lem:addi} and
  \eqref{eq:add3}.
  Since 
  $\overline{F\circ \mu_j}=\sigma(\mu_j)^{-1}F\hat\rho^{}_F(\tilde\mu_j)$ 
  and $\hat\rho^{}_F(\tilde\mu_j)\in SU(2)$, 
  we have
  $f\circ \mu_j=\sigma(\mu_j)^{-1}f(\sigma(\mu_j)^{-1})^*$.
  This implies the last assertion.
\end{pf}

The following lemma plays an important role in later sections:

\begin{lemma}\label{lem:inf}
  Let $F_c\colon{}M\to SL(2,\C)$ be a family of null holomorphic immersions
  such that 
  $\lim_{c\to 0}F_{c}=\id$ and
  $(1/c)F_c^{}F_c^*\in \D_M^{(c)}(G,Q)$.
  Let $l$ be a loop on $M$ and $\tau\in\pi_1(M)$ 
  the deck transformation induced from $l$.
  Then
  \[
    \left.\frac{\partial}{\partial c}\right|_{c=0}
      \rho^{}_{F_{c}}(\tau)=
      \oint_l
          \begin{pmatrix}
             G   & -G^2 \\
             1   & -G\hphantom{^2}
          \end{pmatrix}\frac{Q}{dG} \; \; .
  \]
\end{lemma}
\begin{pf}
  Let $z_0$ be a be point of $l$ in $M$ and 
  $\tilde z_0\in\widetilde M$ a lift of $z_0$.
  Put $\tilde z_1=\tau(\tilde z_0)$.
  Then there exists a lift $\tilde l$ of $l$ joining 
  $\tilde z_0$ and $\tilde z_1$. 
  We set 
  $F'_c=\left(\left. \partial F_{c}/\partial c\right)\right|_{c=0}$.
  By Corollary~\ref{cor:add}, $F_{c}$ is a solution of the equation
  \begin{equation}\label{eq:ode3}
     dF_c=c\, \alpha F_c,\qquad
            \alpha=
              \begin{pmatrix}
                 G &   -G^2 \\
                 1 &   -G\hphantom{^2}
              \end{pmatrix}
                            \,\frac{Q}{dG} \; .
  \end{equation}
  Since $F_0=\lim_{c\to 0}F_{c}=\id$, we have $\rho^{}_{F_0}=\id$,
  and since $dF'_c=\alpha F_c^{}+c\alpha F'_c$, we have  $dF'_0=\alpha$,
  where $'=\partial/\partial c|_{c=0}$.
  Integrating this, we have
  \[
       F'_0(\tilde z_1)=
            \int_{\tilde z_0}^{\tilde z_1}
            \alpha+F'_0(\tilde z_0).
  \]
  Hence
  \begin{equation}\label{eq:diff1}
            \left.\frac{\partial}{\partial c}\right|_{c=0}
                (F_{c}\circ \tau)(\tilde z_0)
                =\oint_l\alpha + 
                        F'_0(\tilde z_0).
  \end{equation}
  On the other hand, 
  \begin{equation}\label{eq:diff2}
     \left.\frac{\partial}{\partial c}\right|_{c=0}
                (F_{c}\circ \tau)(\tilde z_0)
            =
     \left.\frac{\partial}{\partial c}\right|_{c=0}
                ( F_{c}\rho^{}_{F_{c}}(\tau))(\tilde z_0)
            =
     \left(\left.\frac{\partial}{\partial c}\right|_{c=0}
                \rho^{}_{F_{c}}(\tau)\right)
                +F'_0(\tilde z_0)
  \end{equation}          
  holds. 
  By \eqref{eq:diff1} and \eqref{eq:diff2},
  we are done.
\end{pf}
The above lemma is a generalization of \cite[Lemma~3.3]{UY2}
in which the same statement for a loop $l$ surrounding an end is 
shown.

\section{\cmcone\ surfaces of finite total curvature}%
         \label{sec:ftc}

In \cite{UY4}, the hyperbolic correspondence of the Jorge-Meeks $n$-oid 
and \cmcone\ surfaces of genus $0$ with Platonic symmetries  have been
constructed. (See also Section~\ref{sec:long}.)
In this section, we construct \cmcone\ surfaces 
of both genus zero and positive genus in $\hyp$,
which correspond to minimal surfaces in Euclidean $3$-space
$\R^3$ described in \cite{BR}, \cite{Rossman}, \cite{Wohlgemuth}.

In fact, we show that for any complete finite-total-curvature minimal
surface in $\R^3$ satisfying certain conditions, there exists 
a one-parameter family of corresponding \cmcone\ surfaces in $\hyp$.
The conditions are quite general in the sense that 
they are satisfied by almost all known examples of
complete finite-total-curvature minimal surfaces.

Osserman \cite{Osserman} showed that any complete finite-total-curvature
minimal surface in $\R^3$ can be represented as a conformal immersion 
$f\colon{} M \to \R^3$, where $M$ is a compact Riemannian manifold with a
finite
number of points removed, i.e.~$M = \overline{M} \setminus
\{e_1,...,e_r\}$.  
Let 
$S = \{x \in M \; | \; \exists y \neq x ,f(y) = f(x)\}$ 
be the intersection set of $f$.  
The following definition describes minimal surfaces that are
embedded everywhere except in neighborhoods of the ends.  

\begin{definition}\label{def:almost}
  If there exist $r$ disjoint open disks $U_j \subset \overline{M}$
  ($j =1,\dots,r$) such that $e_j \in U_j$ and 
  $S \subset \cup_{j=1}^r U_j$,
  then the immersion $f$ is {\it almost embedded}.
\end{definition}

Next we define what we mean by a symmetric minimal immersion.  

\begin{definition}\label{def:symmetric}
  A minimal immersion $f$ is {\it symmetric\/} if there is a subregion
  $D \subset f(M)$ that is a disk bounded by non-straight planar
  geodesics, and does not contain any non-straight planar
  geodesics in its interior.  
\end{definition}

If $f$ is symmetric with subdisk $D$, then $D$ generates the entire
surface by reflections across planes containing boundary planar
geodesics.  
(Since each curve in the boundary of $D$ is a planar geodesic, 
 the surface can be smoothly extended across these boundary curves, 
 by the Schwartz reflection principle.)  
Note that since $f$ has finite total curvature, $f$ is not
periodic.  Therefore, if any two boundary curves of $D$ lie in
parallel planes, they must actually lie in the same plane.  

\begin{lemma}\label{lem:threeplanes}
  If $f$ is almost embedded, and if $f$ is symmetric with subdisk $D$,
  then there are at most three planes that contain all of the 
  boundary planar geodesics of $D$.  These \rom(at most three\rom) 
  planes lie in general position, and meet pairwise at angles of
  the form ${\pi}/m$, where $m \geq 2$ is an integer.  
\end{lemma}

\begin{pf}
  Let $P_1,\dots,P_n$
  be the smallest set of distinct planes that contain 
  $\partial D$, and let $\mu_j$ ($j=1,\dots,n$) be the reflection
  with respect to the plane $P_j$.  
  Since $f$ is of finite total curvature, 
  the group $\Gamma$ generated by the symmetries
  $\mu_1,\dots,\mu_n$ is finite.
  It is well-known that $\Gamma$ has a fixed point
  and the number of planes of the fundamental chamber is at most three
  (see \cite[Chapters~4--6]{bou}). 
\end{pf}

Let us name these distinct planes $P_1,\dots, P_s$,
where $s$ ($=2$ or $3$) is the number of the planes.
Let the boundary planar geodesics of $D$ contained in $P_j$ 
be called $S_{j,1}$, $S_{j,2}$, \dots, $S_{j,d_j}$  
($j=1,\dots,s$). 
We now define what we mean by
a non-degenerate set of period problems.  
Let $d$ equal the number of smooth boundary planar geodesics 
minus the number of planes.  
Thus $d= d_1+d_2+d_3-3$ if $s=3$, and $d = d_1+d_2-2$ if $s=2$.  

\begin{definition}
  Let $f\colon{}M\to \R^3$ be an almost embedded minimal immersion 
  that is symmetric with subdisk $D$.  
  Then the set of period problems for $f$ is {\it non-degenerate\/} 
  if there exists a continuous $d$-parameter family of disks $D_\lambda$, 
  ($\lambda=(\lambda_1,...,\lambda_{d}),\,|\lambda|<\varepsilon$) such that
  \begin{enumerate}
    \item   $D_{(0,0,...,0)} = D$.
    \item   $\partial D_\lambda = \cup_{j=1}^{s}(\cup_k^{d_j}
              S_{j,k}(\lambda))$  
            such that  each $S_{j,k}(\lambda)$ is 
            a planar geodesic lying in a plane
            $P_{j,k}(\lambda)$ parallel to $P_j$.  
    \item   Letting $\Per_{j,k}(\lambda)$ 
            ($j=1,\dots,s$, $k=2,\dots,d_j$) be the oriented distance 
            between the plane $P_{j,k}(\lambda)$ and $P_{j,1}(\lambda)$, 
            the map 
            $\Lambda\colon{} \lambda = (\lambda_1,\dots,\lambda_{d}) \mapsto
                 (\Per_{j,k}(\lambda)) \in \R^{d}$
            is an open map onto  a neighborhood of $0$.  
  \end{enumerate}
\end{definition}

We reflect $D_\lambda$ infinitely often 
to get a simply connected complete surface
$\widetilde M_\lambda$.  
Let $M$ be the initial minimal surface. 
Then the universal cover $\widetilde M$ of $M$ coincides with
$\widetilde M_0$.
The initial fundamental disk $D_\lambda$ is contained in 
$\widetilde M_\lambda$.
Also we have associated
reflections $\tilde{\mu}_{\lambda,j,k}$ on $\widetilde M_\lambda$
with respect to $D_\lambda$, 
and we  have the properties 
\begin{equation} \label{eq:gr_rel}
  \tilde{\mu}_{\lambda,j,k} \circ \tilde{\mu}_{\lambda,j,k}
      = \id \; , 
  \qquad
  (\tilde{\mu}_{\lambda,j,k} \circ 
  \tilde{\mu}_{\lambda,j^\prime,k^\prime})^{m_{j,k,j^\prime,k^\prime}} = 
  \id \;,
\end{equation}
whenever the planes containing $S_{j,k}$ and $S_{j^\prime,k^\prime}$
meet at an angle of ${\pi}/{m_{j,k,j^\prime,k^\prime}}$.  
Let $\widetilde{\Gamma}(\lambda)$ be 
the group generated by these reflections, 
with the above properties \eqref{eq:gr_rel}.
Since $m_{j,k,j',k'}$ does not depend on $\lambda$, 
$\widetilde{\Gamma}(\lambda)$ is isomorphic to
$\widetilde{\Gamma}(0)$, 
so let 
$\iota\colon{} \widetilde \Gamma(\lambda) \rightarrow  \widetilde{\Gamma}(0)$ 
be the canonical isomorphism.
Let $\Gamma$ be the deck transformation of the 
universal cover $\widetilde M$ of the initial minimal surface $M$.
$\Gamma$ is considered as a subgroup of $\widetilde{\Gamma}(0)$.
Let $\Gamma_\lambda = \iota^{-1}(\Gamma)$, 
and let $M_\lambda = \widetilde{M}_\lambda/{\Gamma_\lambda}$.  
Then $M_\lambda$ is diffeomorphic to 
the initial minimal surface $M=\widetilde M_0/\Gamma$. 
On the initial minimal surface, each reflection 
$\tilde{\mu}_{0,j,k}$ on $\widetilde M$ induces the reflection
${\mu}_{j,k}\colon{}M\to M$. 
Thus for any $\tau \in \Gamma$, 
there exists a $\tau^\prime \in \Gamma$ such that 
$\tilde{\mu}_{0,j,k} \circ \tau = \tau^\prime \circ\tilde{\mu}_{0,j,k}$. 
By the definition of $\Gamma_\lambda$, any $\tau\in \Gamma_\lambda$
has also an element $\tau'\in \Gamma_\lambda$, 
such that 
$\tilde{\mu}_{\lambda,j,k} \circ \tau = \tau^\prime \circ
 \tilde{\mu}_{\lambda,j,k}$. 
This implies that $\tilde{\mu}_{\lambda,j,k}$ also induces
a reflection $\mu_{\lambda,j,k}\colon{}M_\lambda\to M_\lambda$.

Since $D_\lambda$ is a minimal disk bounded by
planar geodesics, 
the Gauss map $G_\lambda$ and the Hopf differential $Q_\lambda$ 
are well-defined on $D_\lambda$. 
Thus $G_\lambda$ and $Q_\lambda$ can be extended on
$\widetilde{M}_\lambda$, 
so that they satisfy the properties
$\overline{G_\lambda \circ \tilde{\mu}_{j,k}} = 
 \sigma(\mu_{j,k})^{-1} * G_\lambda$ 
and $\overline{Q_\lambda \circ \tilde{\mu}_{j,k}} = Q_\lambda$.
Here,
$\sigma({\mu}_{\lambda,j,k}) \in SU(2)$ is explicitly
given by
\begin{equation}\label{eq:rho}
  \sigma({\mu}_{\lambda,j,k})^{-1}
         =\begin{pmatrix}  
            \nu_2 + i\nu_1 & -i\nu_3 \\
            -i\nu_3  & \nu_2-i\nu_1
          \end{pmatrix},
\end{equation}
where $\nu=(\nu_1,\nu_2,\nu_3)$ is  the unit normal vector 
perpendicular to the plane $P_{j,k}$. 
(Since $P_{j,k}$ is parallel to $P_{j,1}$, the unit vector
 $\nu$ does not depend on the choice of $k$.)
One can easily check the above formula as follows:
Since the left hand matrix is in $SU(2)$,
the transformation 
\[
  T\colon{}\C\cup \{\infty\}\ni z
   \mapsto 
   \overline{\frac{(\nu_2 + i\nu_1)z -i\nu_3}{-i\nu_3z +( \nu_2-i\nu_1)}}
   \in \C\cup \{\infty\}
\]
preserves the metric $d\sigma^2_0=4dzd\overline{z}/(1+|z|^2)^2$.
Moreover, the inverse matrix is the conjugate matrix.
Thus $T$ is an  isometric involution.
Let $\widetilde T\colon{}S^2(1)\to S^2(1)$ be the induced
isometric involution.
By straightforward calculation, $\widetilde T$ corresponds to
a reflection with respect to the plane passing through the origin
perpendicular to $\nu$.

\begin{lemma}
  $G_\lambda$ and $Q_\lambda$ are single-valued on $M_\lambda$ for any
  $\lambda$.
\end{lemma}

\begin{pf}
  It suffices to show that 
  $G_\lambda \circ \tau = G_\lambda$ and 
  $Q_\lambda\circ \tau=Q_\lambda$
  for each $\tau\in\Gamma_{\lambda}$.
  Since $\tau$ is an orientation preserving diffeomorphism on 
  $\widetilde M_\lambda$, it can be written as
  $\tau=\tilde{\mu}_{\lambda,j_1,k_1} \circ \dots \circ 
        \tilde{\mu}_{\lambda,j_{2n},k_{2n}}$.
  Since 
  $\overline{Q_\lambda \circ \tilde{\mu}_{\lambda,j,k}} =Q_\lambda$,
  obviously $Q\circ \tau=Q$ holds.
  On the other hand, we have
  \[
      G\circ \tau=
         \left\{
            \sigma({\mu}_{\lambda,j_{1},k_{1}})
            \overline{\sigma({\mu}_{\lambda,j_{2},k_{2}})}
            \cdots
            \sigma({\mu}_{\lambda,j_{2n-1},k_{2n-1}})
            \overline{\sigma({\mu}_{\lambda,j_{2n},k_{2n}})}
         \right\}^{-1}
         * G.
  \]
  Note that $G\circ \tau=G$ holds for $\lambda=0$, 
  because $G$ coincides with the Gauss map of 
  the initial minimal surface when $\lambda=0$.
  This implies that
  \[
      \sigma({\mu}_{\lambda,j_{1},k_{1}})
      \overline{\sigma({\mu}_{\lambda,j_{2},k_{2}})}
      \cdots
      \sigma({\mu}_{\lambda,j_{2n-1},k_{2n-1}})
      \overline{\sigma({\mu}_{\lambda,j_{2n},k_{2n}})}
      =\pm\id
  \]
  when $\lambda=0$. 
  Thus it holds for any $\lambda$ because 
  $\sigma({\mu}_{j,k})$ does not depend on $\lambda$, 
  by \eqref{eq:rho}. 
\end{pf}

Thus we have an {\it abstract\/} Riemann surface
$M_{\lambda}$. Moreover, 
the induced conformal metric of $\widetilde{M}_{\lambda}$ defined by
$ds^2_\lambda:= (1+ |G_\lambda|^2)^2{Q_\lambda}/{dG_\lambda}\cdot
    \overline{\left({Q_\lambda}/{dG_\lambda}\right)}$
is single-valued on $M_{\lambda}$.
Since $\widetilde M_{\lambda}$ is complete and $M_{\lambda}$ consists of
a finite number of congruent copies of the finite-total-curvature disk 
$D_\lambda$, 
$M_{\lambda}$ is complete and has finite total curvature. 
Hence $M_{\lambda}$ can be written as
\[
    M_{\lambda}=\overline M_{\lambda}
                \setminus\{e_1,\dots,e_r\}\;,
\]
where $\overline M_{\lambda}$ is a compact Riemann surface, 
and $e_1, \dots, e_r$ are points corresponding
to the ends of $M_{\lambda}$.
Let $f_\lambda\colon{} \widetilde M_\lambda \to \R^3$ be the immersion of the 
(periodic) minimal surface containing $D_\lambda$.
Then $G_\lambda$ and $\omega_\lambda = {Q_\lambda}/{dG_\lambda}$ 
are the Weierstrass data
of $f_\lambda$.  
Since $D_\lambda$ is continuous in $\lambda$, 
we see that $G_\lambda$, and $Q_\lambda$ are also continuous in 
$\lambda$.  
{\it For the sake of simplicity, we express each reflection
    by $\mu_{j,k}$ instead of $\mu_{\lambda,j,k}$.} 
We can place $D_\lambda$ in $\R^3$ so that the boundary curves
$S_{1,k}$ lie in the plane $\{x_2 = 0\}$, 
the boundary curves $S_{2,k}$ lie in a vertical plane containing 
the $x_3$-axis, 
and the boundary curves $S_{3,k}$ lie in a nonvertical plane.  
Since the angle between $P_1$ and $P_2$ is of the form $\pi/m$ 
(see Lemma~\ref{lem:threeplanes}), we have
\begin{equation}\label{eq:JMsymm2}
    \sigma(\mu_{1,k})=\id,\quad
    \sigma(\mu_{2,k})=
        \begin{pmatrix} 
           e^{\pi i/m} & 0 \\
           0 & e^{-\pi i/m} 
        \end{pmatrix},\quad
    \sigma(\mu_{3,k})=
        \begin{pmatrix} 
           \hphantom{-}p_0   &   q_0 \\ 
           -\bar q_0 & \bar p_0 
        \end{pmatrix},
\end{equation}
where $m$ is an integer, $|p_0|^2 + |q_0|^2 = 1$,
and $q_0 \neq 0$, $\Re(q_0)=0$.

\begin{proposition}\label{prop:JMStep2}
  Suppose $m\ge 2$. 
  Then,  for a sufficiently small $\varepsilon>0$,
  there exists a $(d+1)$-parameter family 
  $\{F_{c,\lambda}\}_{|c|<\varepsilon,\lambda \in
          B_\varepsilon(0)}$,
  \rom($B_\varepsilon(0) \subset \R^{d}$\rom)
  of null holomorphic immersions of $\widetilde{M}_\lambda$
  into $SL(2,\C)$ with the following properties.
  \begin{enumerate}
   \item   \label{item:JM1}
           $f_{c,\lambda}=(1/c)F^{}_{c,\lambda}{F_{c,\lambda}^*}\in
            \D_{M_{\lambda}}^{(c)}(G_{\lambda},Q_{\lambda})$ 
           for each $(c,\lambda)$ \rom($c\neq 0$\rom).
   \item   \label{item:JM2}
           $F_{c,\lambda}$ is smooth in $c$ and continuous in $\lambda$.
   \item   \label{item:JM3}
           $\lim_{c\to 0}F_{c,\lambda}=\id$.
   \item   \label{item:JM4}$\lim_{c\to
0}\hat\rho^{}_{F_{c,\lambda}}(\tilde\mu_{j,k})=
           \sigma(\mu_{j,k})$,
           where $\rho^{}_{F_{c,\lambda}}(\tilde\mu_{j,k})$ 
           is the matrix given in Lemma~\ref{lem:repref}.
   \item   \label{item:JM5}
           $\hat\rho^{}_{F_{c,\lambda}}(\tilde\mu_{1,1})=\id$.
   \item   \label{item:JM6}
           There exists a smooth function $\xi=\xi(c,\lambda)$ such that
           $|\xi|=1$ and 
           \[
                \hat\rho^{}_{F_{c,\lambda}}(\tilde\mu_{2,1})
                =
                \begin{pmatrix}
                  \xi & 0 \\ 0 & \xi^{-1}
                \end{pmatrix}.
           \]
   \item   \label{item:JM7}
           When $s=3$,
           \[
                 \hat\rho^{}_{F_{c,\lambda}}(\tilde\mu_{3,1})
                    =\begin{pmatrix}
                       p_{3,1} &   i\beta\\
                       i\beta  &   \overline{p_{3,1}}
                     \end{pmatrix},                       
           \]
           where $p_{3,1}=p_{3,1}(c,\lambda)$
           is a complex-valued function and $\beta=\beta(c,\lambda)$ is
           a real-valued function.
   \item   \label{item:JM8}
           For each $j= 1,\dots, s$, and each $k \geq 2$, 
           \[
               \hat\rho^{}_{F_{c,\lambda}}(\tilde\mu_{j,k})=
                     \begin{pmatrix}
                        p_{j,k} &   i\gamma_{1,j,k}\\
                        i\gamma_{2,j,k}  &   \overline{p_{j,k}}
                     \end{pmatrix}\; ,
           \]                      
          where $p_{j,k}=p_{j,k}(c,\lambda)$ is a complex-valued
          function, 
          and $\gamma_{l,j,k}=\gamma_{l,j,k}(c,\lambda)$
          \rom($l=1,2$\rom) are real-valued functions.
  \end{enumerate}
\end{proposition}

The following lemma is easily proved in view of the
fact that $a\overline{a}=\id$ implies $\overline{a}=a^{-1}$.

\begin{lemma}\label{lem:newalg}
  Let $a$ be a matrix in $SL(2,\C)$.
  Then $a\overline{a}=\id$ holds if and only if $a$ is of the form
  \[
     a=\begin{pmatrix} 
         p & i\gamma_1 \\
         i\gamma_2 & \overline{p} 
       \end{pmatrix}.
  \]
\end{lemma}

\begin{pf*}{{\it Proof of Proposition}~\ref{prop:JMStep2}}
  We prove the proposition in the case of $s=3$.
  When $s=2$, the assertions are proved
  if we simply ignore the discussions after \eqref{eq:JMmat1} in Step~III.
  \renewcommand{\qed}{}
\end{pf*}
\subsection*{Step I}
  Let $z_{0,\lambda}$ be a fixed point in the curve $S_{1,1}$, 
  which depends continuously on $\lambda$.
  Consider the initial value problem according to Corollary~\ref{cor:add}:
  \begin{equation}\label{eq:odefam}
       d \check F_{c,\lambda} \check F_{c,\lambda}^{-1}
              =c\alpha_{\lambda},
            \quad
       \check F_{c,\lambda}(\tilde z_{0,\lambda})=\id,
      \qquad\text{where }
      \alpha_{\lambda}=\begin{pmatrix}
                    G_{\lambda} &   -G_{\lambda}^2\\
                    1           &   -G_{\lambda}
                  \end{pmatrix}
                  \frac{Q_{\lambda}}{dG_{\lambda}} \; . 
  \end{equation}
  $G_{\lambda}$ is real on the curve $S_{1,1}$, since $S_{1,1}$ is 
  a planar geodesic in a plane perpendicular to the $x_2$-axis.
  $Q$ is real on the curve $S_{1,1}$, since 
  $\overline{Q_{\lambda}\circ \mu_{1,1}} = Q_{\lambda}$.  
  Thus, the solutions $ \check F_{c,\lambda}$ are also real on 
  $S_{1,1}$.
  For the solutions $ \check F_{c,\lambda}$ of \eqref{eq:odefam},
  $ \check f_{c,\lambda}=(1/c) \check F_{c,\lambda}^{} 
    \check F_{c,\lambda}^*$
  satisfies \ref{item:JM1}, \ref{item:JM2}.
            
  For $c=0$, \eqref{eq:odefam} reduces to $d \check F_{0,\lambda}=0$.
  So $ \check F_{0,\lambda}$ satisfies \ref{item:JM3}.
  Moreover, since the $\check F_{c,\lambda}$ are real on $S_{1,1}$, 
  $ \overline{\check F_{c,\lambda}\circ\mu_{1,1}}=
    \check F_{c,\lambda}$ on $\widetilde S_{1,1}$,
  where $\widetilde S_{1,1}$ is the lift of $S_{1,1}$
  that contains the point $\tilde z_0$.
  So by the holomorphicity of $\check F_{c,\lambda}$,
  \[
       \overline{ \check F_{c,\lambda}\circ\mu_{1,1}}=
                  \check F_{c,\lambda}\qquad
       \text{ on $\widetilde M_\lambda$}\;.
  \]
  This shows that $\hat\rho^{}_{\check F}(\tilde\mu_{1,1})$ is the identity
  and  \ref{item:JM5} is proved.
  By  \eqref{eq:odefam}, we have $\check F_{0,\lambda}=\id$.
  This proves \ref{item:JM3}.
  By 
  $\overline{G_\lambda \circ \tilde{\mu}_{j,k}} =
      \sigma({\mu}_{j,k})^{-1} * G_\lambda$,
  we have
  \[
     \overline{
       \check F_{0,\lambda} \circ \tilde{\mu}_{j,k} 
     }
     =\sigma({\mu}_{j,k})^{-1}
     \check F_{0,\lambda}
     \hat\rho^{}_{\check F}({\tilde\mu}_{j,k}),
  \]
  where $\sigma(\mu_{j,k})$ is the matrix given by 
  \eqref{eq:rho}.
  Since $\check F_{0,\lambda}=\id$, we have
  $\lim_{c\to 0}\hat\rho^{}_{\check F_{c,\lambda}}(\tilde\mu_{j,k})
    =\sigma(\mu_{j,k})$.
  This proves \ref{item:JM4}.
\subsection*{Step II}
  By \eqref{lem:id}, we have
  \begin{equation}\label{eq:conj2a}
    {\hat\rho^{}_{\check F_{c,\lambda}}(\tilde\mu_{j,k})}
    \overline{
       \hat\rho^{}_{\check F_{c,\lambda}}(\tilde\mu_{j,k})}=\id \; .
  \end{equation}
  By Lemma~\ref{lem:newalg}, we have the following expression:
  \[
     \hat\rho^{}_{\check F_{c,\lambda}}(\tilde\mu_{2,1})
            =\begin{pmatrix}
                 p(c,\lambda) & i\gamma_1 \\
                 i\gamma_2 & \overline{p(c, \lambda)}
             \end{pmatrix}
              \left(\mathstrut \in SL(2,\C)\,\right ).
  \]
  Since we suppose that $m\ge 2$ and 
  $\lim_{c\to 0} p(c,\lambda)=e^{\pi i/m }$,
  there exists a sufficiently small neighborhood $U\subset \R^k$ and
  a positive $\varepsilon>0$ such that the imaginary
  part $\Im(p(c,\lambda))>0$  and $|\Re(p)|<1$ for any $|c|<\varepsilon$ 
  and $\lambda \in U$.
  In such a small neighborhood,
  $\hat\rho^{}_{ F_{c,\lambda}}(\tilde\mu_{2,1})$ has 
  two distinct eigenvalues  $\xi(c,\lambda)$ and $1/\xi(c,\lambda)$.
  Since $|\Re(p)|<1$, $\xi(c,\lambda)$ is not a real number,
  and thus $|\xi(c,\lambda)|=1$ by straightforward calculation.
  Moreover $\xi=\xi(c,\lambda)$ is a smooth function. 
  Now we define a real matrix $u(c,\lambda)\in SL(2,\R)$ by
  \[
      u(c,\lambda)=\left\{
                      1+\frac{\gamma_1\gamma_2}
                      {\left( \Im(p)+\sqrt{1-\Re(p)^2}\right)^2}
                   \right \}^{-1}
                   \begin{pmatrix}
                     1 & {-i\gamma_1}/(p-\overline{\xi}) \\
                     {-i\gamma_2}/(\overline{p}-\xi) & 1
                   \end{pmatrix}.
  \]
  Since $\Im(p(c,\lambda))>0$, we can easily see that
  $\overline{p}-\xi$ is a non-vanishing imaginary number 
  for any $\lambda\in U$ and $|c|<\varepsilon$.
  Thus $u=u(c,\lambda)$ is defined as a smooth $SL(2,\R)$-valued
  function. 
  By straightforward calculation,
  we have 
  \begin{equation}\label{eq:diag}
    \hat\rho^{}_{\check F_{c,\lambda}}(\tilde\mu_{2,1})
      =u(c,\lambda)
       \begin{pmatrix} 
          \xi(c,\lambda) & 0 \\
          0 & \overline{\xi(c, \lambda)}.
       \end{pmatrix}
       u(c,\lambda)^{-1}.
  \end{equation}
  We set $\hat F_{c,\lambda}=\check F_{c,\lambda}u(c,\lambda)$
  and $\hat f_{c,\lambda}=(1/c)\hat F_{c,\lambda}^{}\hat F_{c,\lambda}^*$.
  Since
  \[
      \overline{
         \check F_{0,\lambda} \circ \tilde{\mu}_{j,k} 
      }
      =\sigma({\mu}_{j,k})^{-1}
       \check F_{0,\lambda}\hat\rho^{}_{\check F}({\tilde\mu}_{j,k})
  \]
  and $u(c,\lambda)\in SL(2,\R)$,
  we have by Corollary~\ref{lem:repchange}
  \[
      \hat\rho^{}_{\hat F_{c,\lambda}}(\mu_{j,k})={u(c,\lambda)}^{-1} 
      \hat\rho^{}_{\check F_{c,\lambda}}(\mu_{j,k})
       \overline{u(c,\lambda)}\qquad (j=1,\dots,s,\,\, k=1,\dots,d_j).
  \]
  By \eqref{eq:diag}, we have \ref{item:JM6} for $\hat f_{c,\lambda}$.
  Since $u(0,\lambda)=\id$, $\hat f_{c,\lambda}$
  satisfies the same properties \ref{item:JM1}--\ref{item:JM5} as 
  $\check f_{c,\lambda}$.
\subsection*{Step III}
  By \eqref{eq:conj2a} and Lemma~\ref{lem:newalg},
  the matrices $\hat\rho^{}_{\hat F_{c,\lambda}}(\tilde\mu_{3,1})$ 
  and $\hat\rho^{}_{\hat F_{c,\lambda}}(\tilde\mu_{j,k})$
  ($j\ge2$) are written in the form
  \begin{equation}\label{eq:JMmat1}
     \hat\rho^{}_{\hat F_{c,\lambda}}(\tilde\mu_{3,1})
        =\begin{pmatrix}
            p_{3,1}   &   i\beta_1 \\
            i\beta_2   & \overline{p_{3,1}}
         \end{pmatrix},
            \qquad \quad
       \hat\rho^{}_{\hat F_{c,\lambda}}(\tilde\mu_{j,k})
        =\begin{pmatrix}
             p_{j,k}   &   i\gamma_{1,j,k} \\
             i\gamma_{2,j,k}    & \overline{p_{j,k}}
         \end{pmatrix} \; ,
  \end{equation}
  where $\beta_1$, $\beta_2$, $\gamma_{1,j,k}$ and $\gamma_{2,j,k}$ 
  are real-valued functions of $c$ and $\lambda$.
  Since 
  $\lim_{c\to 0}\hat\rho^{}_{\hat F_{c,\lambda}}
         (\tilde\mu_{3,1})=\sigma(\mu_{3,1})$,
  we have
  \begin{equation}\label{eq:JMcomp}
     p_{3,1}(0,\lambda)=p_0, \qquad
     \qquad
     \beta_1(0,\lambda)=\beta_2(0,\lambda)=-iq_0 \neq 0 \; . 
  \end{equation}
  By \eqref{eq:JMcomp},
  there exists a positive number $\varepsilon$ such that 
  $\beta_1/\beta_2>0$ holds for $|c|<\varepsilon$.
  For $c$ in such a range, we set 
  ${f_{c,\lambda}}=(1/c)F_{c,\lambda}^{} F_{c,\lambda}^*$, 
  where
  \[
     {F_{c,\lambda}}=
          \hat F_{c,\lambda} 
          \begin{pmatrix}   
            t_{c,\lambda}   &   0 \\
            0   &   t_{c,\lambda}^{-1}
          \end{pmatrix}\;,
          \qquad
         t_{c,\lambda}=\sqrt[4]{\frac{\beta_1}{\beta_2}}
         \,.
  \]
  Since $t_{0,\lambda}=1$ and $t=t_{c,\lambda}$ is a real-valued function,
  $F_{c,\lambda}$ and $f_{c,\lambda}$ satisfy \ref{item:JM1}--\ref{item:JM6}.
  By Corollary~\ref{lem:repchange},
  \[
     \hat\rho^{}_{F_{c,\lambda} }(\tilde\mu_{3,1})=
        \begin{pmatrix}   
          t_{c,\lambda}^{-1}  &   0 \\
          0   &   t_{c,\lambda}
        \end{pmatrix}
        \hat\rho^{}_{\hat F_{c,\lambda}}(\tilde\mu_{3,1})
        \begin{pmatrix}   
          t_{c,\lambda}   &   0 \\
          0   &   t_{c,\lambda}^{-1}
        \end{pmatrix}=
        \begin{pmatrix}
          p_{3,1} &   i\beta \\
          i\beta & \overline{p_{3,1}}
        \end{pmatrix} \; ,
  \]
  where $\beta=\sqrt{\beta_1\beta_2}$.
  Thus $f_{c,\lambda}$  satisfies \ref{item:JM7}.
  By \eqref{eq:JMmat1}, we have
  \[
       \hat\rho^{}_{{F}}(\tilde\mu_{j,k})=
       \begin{pmatrix}   
          t_{c,\lambda}^{-1}  &   0 \\
          0   &   t_{c,\lambda}
       \end{pmatrix} 
       \hat\rho^{}_{\hat F_{c,\lambda}}(\tilde\mu_{j,k})
       \begin{pmatrix}   
          t_{c,\lambda}   &   0 \\
          0   &   t_{c,\lambda}^{-1}
       \end{pmatrix}=
       \begin{pmatrix}
          p_{j,k}   &   it_{c,\lambda}^{-2}\gamma_{1,j,k} \\
          it_{c,\lambda}^{2}\gamma_{2,j,k}    & \overline{p_{j,k}}
       \end{pmatrix}. 
  \] 
  Replacing $\gamma_{1,j,k}$ and $\gamma_{2,j,k}$ by
  $t_{c,\lambda}^{-2}\gamma_{1,j,k}$ and $t_{c,\lambda}^2\gamma_{2,j,k}$,
  we get \ref{item:JM8}.
\qed

\begin{remark}\label{rem:converge}
  In the Poincar\'e model of radius $1/c$, 
  the immersion $f_{c,\lambda}$ converges to the initial
  minimal immersion $f_\lambda$. 
  Indeed, by \eqref{eq:odefam} we have  
  $d F_{c,\lambda} F_{c,\lambda}^{-1}=
     d \check F_{c,\lambda} \check F_{c,\lambda}^{-1}=c\alpha_\lambda$.
  By the same argument as in Lemma \ref{lem:inf},
  we have $d F'_{0,\lambda}=\alpha_{\lambda}$,
  where
  $F'_{0,\lambda}=\left(\partial F_c/\partial c\right)|_{c=0}$.
  In the complexified Poincar\'e plane of radius $1/c$,
  $F_{c,\lambda}$ converges to $F'_{0,\lambda}$ by Lemma 2.1 of \cite{UY2}.
  On the other hand,
  the initial minimal immersion $f_\lambda$ can be
  expressed as $f_\lambda=F'_{0,\lambda}+\overline{F'_{0,\lambda}}$,
  where we identify a point $(x_1,x_2,x_3)$ in
  $\C^3$ with the matrix 
  \[
      \begin{pmatrix}
           x_3 & x_1+ix_2 \\
           x_1-ix_2 & -x_3
      \end{pmatrix}
  \]
  in the Lie algebra of $SL(2,\C)$.
  Thus by Corollary 2.2 of \cite{UY2}, we can conclude
  $f_{c,\lambda}\to f_{\lambda}$.
\end{remark}

It should be remarked that
$\hat\rho^{}_{F_{c,\lambda}}(\tilde\mu_{1,1})$,
$\hat\rho^{}_{F_{c,\lambda}}(\tilde\mu_{2,1})$,
$\hat\rho^{}_{F_{c,\lambda}}(\tilde\mu_{3,1}) \in SU(2)$.
Moreover,
$\hat\rho^{}_{F_{c,\lambda}}(\tilde\mu_{j,k}) \in SU(2)$
if and only if $\gamma_{1,j,k} = \gamma_{2,j,k}$.
By Proposition~\ref{prop:su2cond},
the period problem in the hyperbolic case reduces to showing that
$\gamma_{1,j,k} = \gamma_{2,j,k}$.  
To show that this can be done, 
we have the following lemma.  Recall that $o(c)$ denotes any function of
$c$ that tends to zero faster than $c$ itself as $c \rightarrow 0$.  

\begin{lemma}\label{lem:JMinter2}
  Let $F_{c,\lambda}$ be as in Proposition~\ref{prop:JMStep2}. 
  Then
  \[
       (\gamma_{1,j,k}-\gamma_{2,j,k})(c,\lambda) = 
            2c\cdot \Per_{j,k}(\lambda) + o(c)
        \qquad
        (j=1,\dots,s\,,\,\,k=2,\dots,d_j).
  \]
\end{lemma}
\begin{pf}
  Let $l_{j,k}^{(1)}$ be a curve in $D_\lambda$ 
  starting from a point  on $S_{j,1}$ and 
  ending at a point on $S_{j,k}$, and let
  $l_{j,k}^{(2)}$ be  its reflected curve
  across the plane containing $S_{j,1}$.  
  Then the curve $l_{j,k}$ obtained as the composite of
  the reversed oriented
  curve $-l_{j,k}^{(2)}$ with $l_{j,k}^{(1)}$ can be
  considered as a closed loop in $M_\lambda$.
  Let $\tau_{j,k}$ be the 
  corresponding  element of $\pi_1(M_\lambda)$.
  It can be identified with an element of the deck
  transformations of $\widetilde M_\lambda$, 
  if we choose the base point in $\widetilde M_\lambda$
  as the initial point of $l_{j,k}$.
  Then we have
  $\tau_{j,k}= \tilde\mu_{j,k}\circ\tilde\mu_{j,1}$.  
  We have from Remark~\ref{rem:refrel} and \eqref{eq:Add} that 
  \[
     F_{c,\lambda}\circ\tau_{j,k}
       =F_{c,\lambda}\circ\tilde\mu_{j,k}\circ\tilde\mu_{j,1}
       =F_{c,\lambda}
        {\hat\rho^{}_{F_{c,\lambda}}(\tilde\mu_{j,k})}
                   \overline{ \hat\rho^{}_{F_{c,\lambda}}(\tilde\mu_{j,1})}.
  \]
  By the Weierstrass representation
  formula for the initial minimal surface with
  the Weierstrass data $(G_\lambda,\omega_\lambda)$,
  (where $\omega_\lambda=Q_\lambda/dG_\lambda$),  
  we have
  \begin{equation}\label{eq:JMint}
     \Re\oint_{l_{j,k}}
     \left(
        1-G_{\lambda}^2, i(1+G_{\lambda}^2),
        2G_{\lambda}
     \right)\omega_{\lambda}  
       =  2\cdot \Per_{j,k}(\lambda)  \cdot 
                    (\nu_1,\nu_2,\nu_3)  \, ,
  \end{equation}
  where $(\nu_1,\nu_2,\nu_3)$ is a unit vector perpendicular to the
  planes containing $S_{j,1}$ and $S_{j,k}$.  
  If we set 
  \[
       \Im\oint_{l_{j,k}}
         \left(
            1-G_{\lambda}^2, i(1+G_{\lambda}^2),
            2G_{\lambda}
         \right)\omega_{\lambda}  
         =  (\eta_1,\eta_2,\eta_3) \, ,
  \]
  then we have from Lemma~\ref{lem:inf}
  \begin{equation}\label{eq:star}
     \left.\frac{\partial}{\partial c}\right|_{c=0}
           \rho^{}_{F_{c,\lambda}}(\tau_{3,k})=
     \oint_{l_{3,k}}
          \begin{pmatrix}
              G_{\lambda} &   -G_{\lambda}^2\\
              1       &   -G_{\lambda}
          \end{pmatrix}
          \omega_{\lambda} 
     = \Per_{j,k}(\lambda)
        A + iB,
  \end{equation}
  where $A$ and $B$ are the matrices given by
  \begin{equation}\label{eq:AB}
     A=\begin{pmatrix} 
         \nu_3 & \nu_1+i\nu_2 \\
         \nu_1-i\nu_2& -\nu_3 
       \end{pmatrix}, \qquad
     B=\begin{pmatrix} 
         \eta_3 & \eta_1+i\eta_2 \\
         \eta_1-i\eta_2& -\eta_3 
       \end{pmatrix}\,.
  \end{equation}
  Since 
  $\rho^{}_{F_{c,\lambda}}(\tau)=
    \hat\rho^{}_{F_{c,\lambda}}(\tilde\mu_{j,k})
    \overline{\hat\rho^{}_{F_{c,\lambda}}(\tilde\mu_{j,1})}$, 
  by \ref{item:JM4} of Proposition~\ref{prop:JMStep2},
  we have
  \[
       \rho^{}_{F_{c,\lambda}}(\tau_{3,k})'=
      {(\rho^{}_k)'} \overline{\sigma} + {\sigma}
      \overline{(\rho^{}_1)'},
  \]
  where we simplify the notation as  
  $\rho^{}_k=\hat\rho^{}_{F_{c,\lambda}}(\tilde\mu_{j,k})$,
  $\sigma=\sigma(\mu_{j,1})=\sigma(\mu_{j,k})$
  and the prime denotes  differentiation with respect to $c$
  evaluated at $c=0$. 
  Since ${\rho^{}_j}\overline{\rho^{}_j}=\id$, 
  we have
  ${(\rho^{}_j)}'\overline{\sigma}+{\sigma}\overline{(\rho^{}_j)'}=0$.
  Thus by \eqref{eq:star} we have
  $\rho^{}_{F_{c,\lambda}}(\tau_{3,k})'=
   \{(\rho^{}_k)'-(\rho^{}_1)'\}\overline{\sigma}$.
  By Lemma~\ref{lem:id}, $\sigma=\overline{\sigma}^{-1}$ holds.
  Thus we have
  \[
      (\rho^{}_k)'-(\rho^{}_1)'=\Per_{j,k}(\lambda) A\sigma + i B\sigma.
  \]
  For a matrix $a=(a_{ij})_{i,j=1,2}$, we define a function
  $\Delta[a]$ by $\Delta[a]=a_{12}-a_{21}$.
  Since $\Delta[\rho^{}_1]=0$, we have $\Delta[(\rho^{}_1)']=0$.
  Hence
  \[
      i\{(\gamma_{1,j,k})'-(\gamma_{2,j,k})'\} =
      \Delta[(\rho^{}_k)']
      =\Delta[(\rho^{}_k)']-\Delta[(\rho^{}_1)'] 
      =\Per_{j,k}(\lambda)\Delta[A\sigma] 
      + i\Delta[B\sigma].
  \]
  By \eqref{eq:rho} and \eqref{eq:AB}, we can directly compute 
  that $\Delta[A\sigma]=2i$ and 
  $\Delta[B\sigma]\in i\R$.  
  Since $(\gamma_{1,j,k})'-(\gamma_{2,j,k})'\in \R$,
   we have
  \begin{equation}\label{eq:a}
      (\gamma_{1,j,k})'-(\gamma_{2,j,k})'=2\cdot
      \Per_{j,k}(\lambda).
  \end{equation}
  On the other hand, by \ref{item:JM4} of
  Proposition~\ref{prop:JMStep2}, 
  it holds that
  \begin{equation}\label{eq:b}
      \lim_{c\to 0}(\gamma_{1,j,k}-\gamma_{2,j,k})
      =\lim_{c\to 0}\Delta(\rho^{}_j)=\Delta(\sigma)=0.
  \end{equation}
  By \eqref{eq:a} and \eqref{eq:b},
  we have the conclusion.
\end{pf}

Thus we have equated the period problems for the minimal surface in
$\R^3$ with the period problems for the \cmcc{} surface in
$\hypc$.  If the period problem on the minimal surface 
is non-degenerate, 
the period map $\lambda\mapsto \Per_{j,i}(\lambda)$ 
is an open map onto a neighborhood of the origin.  
Since $F$ is smooth in $c$ and continuous in $\lambda$, 
for any sufficiently small $c$ there exists a $\lambda_c$ such that 
$\gamma_{2,j,k}(c,\lambda_c)=\gamma_{1,j,k}(c,\lambda_c)$ for
all $k$ and $j$.  
For such a pair $(c,\lambda_c)$, 
we see that
$f_{c,\lambda_c}$ satisfies Proposition~\ref{prop:su2cond}.
Thus the \cmcc{} surface $f_{c,\lambda_c}$
is an element of $I_{M}^{(c)}(G_{\lambda_c},Q_{\lambda_c})$.
By Lemma~\ref{lem:Dempty},
the surface has a complete induced metric, whenever the Hopf differential
$Q$ has  a pole of order at most $2$. 
Moreover, in this case the surface has finite total curvature, 
because its secondary Gauss map has at most a pole at each end.
Multiplying the Poincar\'e model by  $c$, $\hypc$ becomes
$\hyp$ and $f_{c,\lambda_c}$ becomes a \cmcone{} surface.  
As we can do this for any $c$
sufficiently close to 0, we have found a one-parameter family of
\cmcone{} surfaces in $\hyp$.  
Thus we get the following:
\begin{theorem}\label{mainthm}
  Let $f_0$ be a complete almost-embedded 
  minimal surface in $\R^3$ with 
  finite total curvature and non-degenerate period problem,
  and suppose that $f_0$ is symmetric with respect to
  a subdisk $D$.
  Then there
  exists a one-parameter family of 
  {\em \cmcone} surfaces $f_c$
  in $\hyp$ with the same reflectional symmetry properties.
  Moreover it has the following properties\rom:
  \begin{enumerate}
   \item If the Hopf differential $Q$ of the initial minimal
         surface has poles of order at most $2$, then 
         the surfaces also have complete induced metrics with
         finite total curvature.
   \item The \cmcone{} surfaces are irreducible if $s=3$, where
         $s$ is the number of 
         planes containing the boundary of the minimal subdisk $D$.
  \end{enumerate}
\end{theorem}

The final statement follows as follows:
Suppose that one of the \cmcone\ surfaces are reducible.
Then by definition, $\hat\rho_{F_{c,\lambda}}(\tilde \mu_{2,1})$
and $\hat\rho_{F_{c,\lambda}}(\tilde \mu_{3,1})$ must be
commutative,
a contradiction to the fact that
$\hat\rho_{F_{c,\lambda}}(\tilde \mu_{3,1})$ is not diagonal.

Now we apply this theorem to construct several examples.  

\subsection{Jorge-Meeks $n$-oid with one handle}\label{sub:JM}
  Let $n\geq 3$ be an integer.
  Then there exist complete minimal surfaces of genus one
  with $n$ catenoid ends in $\R^3$ (cf.\ \cite{BR}).
  These surfaces look like Jorge-Meeks surface \cite{JM} with
  one attached handle, and so they admit $D_n\times \Z_2$ symmetry,
  where $D_n$ is the dihedral group of order $n$.
  The fundamental piece of the surface is 
  given in Figure~\ref{fig:JMpiece}.    
  As shown in Berglund-Rossman~\cite{BR}, the period problem of this example
  is non-degenerate.
  (Here  $s=3$, $d=1$ and so the period problem is one-dimensional.)
  Thus,
  by Theorem \ref{mainthm}, we have the existence of a one-parameter family 
  of \cmcone{} irreducible $n$-oid cousins in $\hyp$ with 
  finite total curvature and symmetry $D_n\times \Z_2$.

\begin{figure}
\begin{center}
\unitlength=0.6pt
\begin{picture}(454.00,173.00)(56.00,632.00)
\put(150.00,660.00){\vector(0,1){40.00}}
\put(150.00,660.00){\vector(-1,0){40.00}}
\put(150.00,660.00){\vector(-3,-4){15.00}}
\put(286.00,750.00){\vector(1,0){94.00}}
\put(149.00,714.00){\makebox(0,0)[cc]{\footnotesize $x_1$}}
\put(96.00,660.00){\makebox(0,0)[cc]{\footnotesize $x_2$}}
\put(127.00,632.00){\makebox(0,0)[cc]{\footnotesize $x_3$}}
\put(225.00,719.00){\makebox(0,0)[cc]{$D_\lambda$}}
\put(265.00,718.00){\makebox(0,0)[cc]{\footnotesize $S_{1,1}$}}
\put(249.00,652.00){\makebox(0,0)[cc]{\footnotesize $S_{3,2}$}}
\put(190.00,642.00){\makebox(0,0)[cc]{\footnotesize $S_{2,1}$}}
\put(181.00,780.00){\makebox(0,0)[cc]{\footnotesize $S_{3,1}$}}
\put(331.00,729.00){\makebox(0,0)[cc]{\footnotesize reflection}}
\put(332.00,712.00){\makebox(0,0)[cc]{\footnotesize across $S_{3,1}$}}
\put(424.00,718.00){\makebox(0,0)[cc]{\footnotesize $S_{3,1}$}}
\put(421.00,641.00){\makebox(0,0)[cc]{\footnotesize $S_{2,1}$}}
\put(479.00,655.00){\makebox(0,0)[cc]{\footnotesize $S_{3,2}$}}
\put(466.00,725.00){\makebox(0,0)[cc]{\footnotesize $S_{1,1}$}}
\bezier136(236.00,752.00)(249.00,793.00)(250.00,700.00)
\bezier80(250.00,700.00)(254.00,645.00)(260.00,670.00)
\bezier51(260.00,670.00)(230.00,671.00)(230.00,650.00)
\bezier85(230.00,650.00)(183.00,643.00)(190.00,680.00)
\bezier192(190.00,680.00)(235.00,766.00)(140.00,780.00)
\bezier3(137.00,774.00)(137.00,773.00)(135.00,771.00)
\bezier5(135.00,765.00)(136.00,763.00)(139.00,761.00)
\bezier8(147.00,757.00)(152.00,754.00)(155.00,755.00)
\bezier7(166.00,753.00)(171.00,752.00)(173.00,752.00)
\bezier10(183.00,750.00)(189.00,749.00)(193.00,750.00)
\bezier9(204.00,749.00)(209.00,748.00)(213.00,749.00)
\bezier5(224.00,749.00)(228.00,749.00)(229.00,749.00)
\bezier113(471.00,752.00)(481.00,787.00)(480.00,710.00)
\bezier93(480.00,710.00)(482.00,644.00)(490.00,670.00)
\bezier55(490.00,670.00)(457.00,672.00)(460.00,650.00)
\bezier91(460.00,650.00)(421.00,632.00)(420.00,680.00)
\bezier179(420.00,680.00)(474.00,759.00)(390.00,770.00)
\bezier57(495.00,685.00)(459.00,681.00)(465.00,660.00)
\bezier30(465.00,660.00)(473.00,674.00)(460.00,684.00)
\bezier52(460.00,684.00)(440.00,697.00)(420.00,680.00)
\bezier47(495.00,685.00)(503.00,696.00)(500.00,730.00)
\bezier75(500.00,730.00)(498.00,789.00)(510.00,800.00)
\bezier10(494.00,803.00)(488.00,803.00)(484.00,805.00)
\bezier9(464.00,804.00)(458.00,804.00)(455.00,803.00)
\bezier10(434.00,798.00)(429.00,797.00)(424.00,795.00)
\bezier8(409.00,788.00)(404.00,786.00)(401.00,784.00)
\bezier6(395.00,777.00)(391.00,775.00)(391.00,773.00)
\bezier6(389.00,765.00)(389.00,762.00)(389.00,759.00)
\bezier10(395.00,753.00)(399.00,749.00)(404.00,750.00)
\bezier7(415.00,747.00)(419.00,746.00)(422.00,746.00)
\bezier8(432.00,746.00)(436.00,746.00)(440.00,746.00)
\bezier7(451.00,745.00)(455.00,745.00)(458.00,745.00)
\bezier3(463.00,745.00)(465.00,745.00)(466.00,745.00)
\end{picture}
\end{center}
    \caption{The fundamental piece $D_{\lambda}$.}%
    \label{fig:JMpiece}
\end{figure}
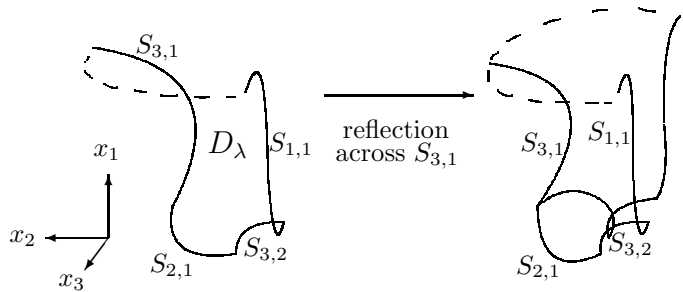

\begin{figure}
\vspace*{-0.5in}
\vspace{2in}
\caption{A genus-1 CMC-1 trinoid in $H^3(-1)$.}
\centerline{\footnotesize
(This picture is due to 
 made by Katsunori Sato of Tokyo Institute of Technology.)}
\end{figure}

\subsection{Genus-0 Jorge-Meeks $n$-oid}

  The genus zero version of the previous example can be shown to exist
  as a \cmcone{} surface in $\hyp$ by the same method.  
  It has been constructed in \cite{UY4}.  
  (There is no period problem ($s=3$, $d=0$), 
   so the non-degenerate condition is trivially satisfied.)
  We examine this example in more detail in the next section.  

\subsection{Genus-0 and higher genus Platonoids}\label{sub:Platon}

  With only slight modifications of the two previous examples, 
  we have \cmcc\ surfaces with Platonic symmetries, which
  correspond to the minimal surfaces in 
  Section 4 of \cite{BR}.
  Higher genus Platonoids were also constructed by 
  Berglund-Rossman \cite{BR}.
  Since they have non-degenerate period problem,
  we get the following irreducible examples by
  Theorem \ref{mainthm}:

  \begin{itemize}
    \item   a \cmcc\ surface of genus $3$ with $4$-ends and 
            tetrahedral symmetry,
    \item   a \cmcc\ surface of genus $5$ with $8$-ends and 
            octahedral symmetry,
    \item   a \cmcc\ surface of genus $7$ with $6$-ends and 
            octahedral symmetry,
    \item   a \cmcc\ surface of genus $11$ with $20$-ends and 
            icosahedral symmetry, and 
    \item   a \cmcc\ surface of genus $19$ with $20$-ends and 
            icosahedral symmetry.  
  \end{itemize}

\subsection{Enneper's surface}

  Enneper's cousin $f$ was found in \cite{Bryant},
  which is a complete simply connected \cmcone{} immersion
  on $\C$ with the same Weierstrass data $(z,dz)$ as 
  the original Enneper surface.
  Since it has no period problem, we can trivially apply
  Theorem~\ref{mainthm} and get the dual $f^\#$
  of Enneper's cousin.
  Since the hyperbolic Gauss map of $f$ is given by
  $G=\tanh z$
  (see \cite{Bryant}), 
  the dual surface $f^\#$ is complete.
  (Since $G$ has an essential singularity at infinity,
   the total curvature of $f^\#$ is infinite.)
  Compared to Enneper's cousin, 
  its dual surface $f^\#$ is aesthetically more appealing
  (see Figure~\ref{fig:enneper}).
  Note that for
  every planar geodesic on the minimal Enneper's surface, there is a
  corresponding planar geodesic on Enneper's cousin.  
  However, for the
  straight lines in the minimal Enneper's surface, 
  there are not corresponding geodesic lines in Enneper's cousin.  
  This is because our construction utilizes 
  the planar geodesics of the minimal surface, 
  but does not concern itself with straight lines in the
  minimal surface.  
  The result is that the symmetry group of the 
  dual of Enneper's cousin is a subgroup of that of 
  minimal Enneper's surface.

Recently, Sato \cite{Sato} has shown the
  existence of higher genus Enneper-type minimal surfaces for any
  positive genus, where he showed a certain non-degeneracy of
the period problems.  
However, his work does not imply that
the period problems are non-degenerate in our sense,
since his fundamental domain is smaller than ours.

\begin{figure}
\vspace{2in}
    \caption{Half of the dual of Enneper's cousin}
    \label{fig:enneper}
\end{figure}

%
%

%
\begin{figure}
\vspace{2in}
    \caption{A 6-ended genus-2 prismoid in Euclidean space}
    \label{fig:prism}
\end{figure}

\subsection{Higher genus prismoid}\label{sub:prism}
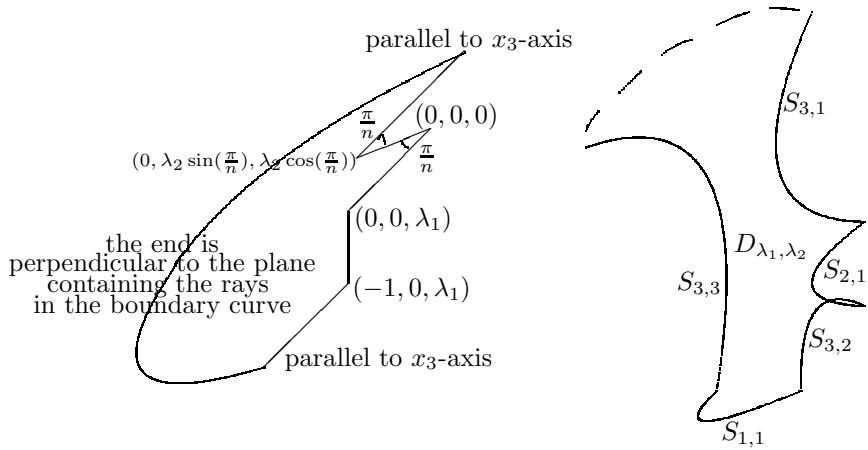
\begin{figure}
\begin{center}
\unitlength=0.8pt
\begin{picture}(378.00,210.00)(42.00,591.00)
\put(180.00,730.00){\line(1,1){51.00}}
\put(180.00,730.00){\line(5,2){35.00}}
\put(215.00,744.00){\line(-1,-1){39.00}}
\put(176.00,705.00){\line(0,-1){34.00}}
\put(176.00,671.00){\line(-1,-1){40.00}}
\put(213.00,730.00){\makebox(0,0)[cc]{\footnotesize $\frac{\pi}{n}$}}
\put(185.00,746.00){\makebox(0,0)[cc]{\footnotesize $\frac{\pi}{n}$}}
\put(227.00,750.00){\makebox(0,0)[cc]{\footnotesize $(0,0,0)$}}
\put(127.00,728.00){\makebox(0,0)[cc]{\tiny $(0,\lambda_2\sin(\frac{\pi}{n}),
\lambda_2\cos(\frac{\pi}{n}))$}}
\put(202.00,701.00){\makebox(0,0)[cc]{\footnotesize $(0,0,\lambda_1)$}}
\put(206.00,668.00){\makebox(0,0)[cc]{\footnotesize $(-1,0,\lambda_1)$}}
\put(195.00,635.00){\makebox(0,0)[cc]{\footnotesize parallel to $x_3$-axis}}
\put(233.00,786.00){\makebox(0,0)[cc]{\footnotesize parallel to $x_3$-axis}}
\put(88.00,690.00){\makebox(0,0)[cc]{\footnotesize the end is}}
\put(88.00,680.00){\makebox(0,0)[cc]{\footnotesize perpendicular to the
plane}}
\put(88.00,670.00){\makebox(0,0)[cc]{\footnotesize containing the rays }}
\put(88.00,660.00){\makebox(0,0)[cc]{\footnotesize in the boundary curve}}
\put(376.00,688.00){\makebox(0,0)[cc]{\footnotesize
$D_{\lambda_1,\lambda_2}$}}
\put(392.00,755.00){\makebox(0,0)[cc]{\footnotesize $S_{3,1}$}}
\put(411.00,676.00){\makebox(0,0)[cc]{\footnotesize $S_{2,1}$}}
\put(342.00,670.00){\makebox(0,0)[cc]{\footnotesize $S_{3,3}$}}
\put(404.00,643.00){\makebox(0,0)[cc]{\footnotesize $S_{3,2}$}}
\put(363.00,599.00){\makebox(0,0)[cc]{\footnotesize $S_{1,1}$}}
\bezier174(231.00,780.00)(134.00,735.00)(94.00,680.00)
\bezier188(94.00,679.00)(42.00,603.00)(136.00,631.00)
\bezier177(395.00,799.00)(351.00,700.00)(420.00,700.00)
\bezier112(420.00,700.00)(370.00,662.00)(420.00,660.00)
\bezier90(420.00,660.00)(388.00,675.00)(390.00,620.00)
\bezier116(390.00,620.00)(320.00,591.00)(350.00,620.00)
\bezier233(350.00,620.00)(372.00,764.00)(288.00,735.00)
\bezier10(288.00,743.00)(290.00,748.00)(294.00,752.00)
\bezier13(304.00,761.00)(311.00,768.00)(314.00,771.00)
\bezier16(324.00,779.00)(332.00,785.00)(338.00,786.00)
\bezier15(350.00,793.00)(358.00,797.00)(365.00,799.00)
\bezier14(375.00,801.00)(383.00,801.00)(389.00,800.00)
\bezier13(190.00,739.00)(192.00,744.00)(193.00,737.00)
\bezier5(201.00,738.00)(202.00,735.00)(204.00,735.00)
\end{picture}
\end{center}
    \caption{The contour $\hat C_{\lambda_1,\lambda_2}$ and the fundamental 
            piece $D_{\lambda_1,\lambda_2}$ for the prismoids}
    \label{fig:prismcont}
\end{figure}

    Let $n\geq 3$ be an integer.  
    There are so called ``higher genus prismoid" minimal surfaces as shown in 
    Figure~\ref{fig:prism} (cf.~\cite{Rossman}).
    These surfaces have $D_n\times \Z_2$ symmetry.
    In these examples, we have two period-killing problems ($s=3$, $d=2$), 
    which correspond to the two parameters $\lambda_1$ and $\lambda_2$ in 
    Figure~\ref{fig:prismcont}.
    The fundamental domain $D_{\lambda_1,\lambda_2}$ 
    as in Figure~\ref{fig:prismcont}
    is a minimal surface in $\R^3$ which is obtained as the conjugate
    surface of the minimal disk bounded by $\hat C_{\lambda_1,\lambda_2}$ in
    Figure~\ref{fig:prismcont}.
    Each segment of the boundary $D_{\lambda_1,\lambda_2}$ 
    is a planar geodesic contained in the plane $P_{j,k}$, where
    $P_{3,1}$, $P_{3,2}$ and $P_{3,3}$ are parallel to the 
    $x_1x_3$-plane, 
    $P_{2,1}$ is perpendicular to
            $(\sin({\pi}/{n}),\cos({\pi}/{n}),0)$, and 
    $P_{1,1}$ is parallel to the $x_1x_2$-plane.
    We define $c_k$ ($k=1,2,3$) by $P_{3,k}=\{x_2=c_k\}$.  
    Then the $c_k$'s are continuous functions in 
    $\lambda_1$ and $\lambda_2$.
    Consider two functions $f_1=c_1-c_3$ and $f_2=c_2-c_3$.
    The following statement implies that the period problem is non-degenerate:
    {\em There exists an open disk in the $(\lambda_1,\lambda_2)$-plane 
    such that its image under the map}
    $(\lambda_1,\lambda_2)
       \mapsto (f_1(\lambda_1,\lambda_2),f_2(\lambda_1,\lambda_2))$
    {\em contains the origin in the $(f_1,f_2)$-plane as an interior point.} 
    This can be proved by arguments found in \cite{BR}. 
    By Theorem \ref{mainthm}, 
    there exist \cmcone{} higher genus prismoid cousins. 

\subsection{Genus-0 prismoids}

  This case is a simpler version of the above example.  In this case,
  there is still a period problem, but now it is only one-dimensional
  ($s=3$, $d=1$), and
  it is known to be non-degenerate (cf.\ \cite{Rossman}).

\subsection{$4$-ended Costa surfaces}\label{sub:costa}
  We refer to a complete minimal surface of 
  finite total curvature 
  with four parallel ends as a ``$4$-ended Costa surface".
  Wohlgemuth~\cite{Wohlgemuth} has constructed several types of
  $4$-ended Costa surfaces, which can have arbitrarily high genus.  
  He has shown that the period problems 
  are non-degenerate for all of his examples.  
  Thus our method applies 
  to all of these examples of Wohlgemuth, 
  producing corresponding \cmcone{} surfaces in $\hyp$.  

  We describe one of the types here 
  (named CSSCFF and CSSCCC by Wohlgemuth).  
  It is a one-parameter family of embedded minimal surfaces of 
  genus $2k-2$ (for any $k \geq 2$), and 
  with four parallel ends.  
  The two outermost ends are always catenoid ends, 
  and the two innermost ends are either catenoid ends 
  (CSSCCC) or planar ends (CSSCFF).  
    
\begin{figure}
\begin{center}
\unitlength=0.8pt
\begin{picture}(300.00,220.00)(90.00,560.00)
\put(160.00,780.00){\line(0,-1){160.00}}
\put(160.00,620.00){\line(-1,0){0.00}}
\put(90.00,740.00){\line(0,-1){160.00}}
\put(90.00,580.00){\line(5,3){70.00}}
\put(160.00,622.00){\line(1,-2){31.00}}
\put(191.00,560.00){\line(0,1){160.00}}
\put(191.00,720.00){\line(-1,2){30.00}}
\put(161.00,780.00){\line(-5,-3){71.00}}
\put(360.00,780.00){\line(0,-1){160.00}}
\put(360.00,620.00){\line(1,-2){30.00}}
\put(390.00,560.00){\line(0,1){160.00}}
\put(390.00,720.00){\line(-1,2){30.00}}
\put(360.00,780.00){\line(-5,-3){70.00}}
\put(290.00,738.00){\line(0,-1){158.00}}
\put(290.00,580.00){\line(5,3){70.00}}
\put(141.00,669.00){\makebox(0,0)[cc]{\footnotesize
$D_{\lambda_1,\lambda_2}$}}
\put(340.00,662.00){\makebox(0,0)[cc]{\footnotesize
$D_{\lambda_1,\lambda_2}$}}
\put(177.00,712.00){\makebox(0,0)[cc]{\footnotesize $S_{1,1}$}}
\put(188.00,615.00){\makebox(0,0)[cc]{\footnotesize $S_{1,2}$}}
\put(118.00,723.00){\makebox(0,0)[cc]{\footnotesize $S_{2,1}$}}
\put(127.00,620.00){\makebox(0,0)[cc]{\footnotesize $S_{2,2}$}}
\put(378.00,714.00){\makebox(0,0)[cc]{\footnotesize $S_{2,1}$}}
\put(390.00,622.00){\makebox(0,0)[cc]{\footnotesize $S_{2,2}$}}
\put(313.00,715.00){\makebox(0,0)[cc]{\footnotesize $S_{3,1}$}}
\put(310.00,621.00){\makebox(0,0)[cc]{\footnotesize $S_{3,2}$}}
\put(345.00,576.00){\makebox(0,0)[cc]{\footnotesize $S_{1,1}$}}
\bezier166(390.00,688.00)(353.00,741.00)(390.00,645.00)
\bezier50(290.00,708.00)(313.00,716.00)(320.00,690.00)
\bezier56(320.00,690.00)(325.00,669.00)(360.00,675.00)
\bezier96(360.00,675.00)(382.00,639.00)(377.00,585.00)
\bezier67(377.00,585.00)(346.00,573.00)(325.00,600.00)
\bezier107(325.00,600.00)(330.00,667.00)(290.00,665.00)
\bezier9(385.00,687.00)(380.00,683.00)(377.00,685.00)
\bezier16(368.00,685.00)(360.00,681.00)(352.00,684.00)
\bezier14(338.00,685.00)(329.00,685.00)(324.00,687.00)
\bezier12(307.00,694.00)(302.00,695.00)(297.00,700.00)
\bezier12(381.00,641.00)(374.00,639.00)(369.00,641.00)
\bezier15(352.00,641.00)(343.00,640.00)(337.00,643.00)
\bezier24(322.00,645.00)(303.00,649.00)(302.00,654.00)
\bezier152(190.00,694.00)(156.00,732.00)(190.00,635.00)
\bezier189(90.00,661.00)(175.00,666.00)(90.00,605.00)
\bezier53(160.00,685.00)(173.00,662.00)(174.00,635.00)
\bezier57(174.00,635.00)(172.00,616.00)(191.00,582.00)
\bezier61(161.00,685.00)(134.00,675.00)(122.00,705.00)
\bezier45(122.00,705.00)(121.00,719.00)(90.00,719.00)
\bezier14(183.00,688.00)(175.00,686.00)(169.00,686.00)
\bezier12(152.00,687.00)(144.00,688.00)(140.00,689.00)
\bezier15(124.00,693.00)(116.00,694.00)(110.00,699.00)
\bezier7(102.00,704.00)(99.00,706.00)(97.00,710.00)
\bezier17(183.00,634.00)(175.00,631.00)(166.00,633.00)
\bezier15(151.00,635.00)(142.00,635.00)(137.00,639.00)
\bezier12(122.00,642.00)(117.00,642.00)(111.00,646.00)
\bezier7(103.00,649.00)(99.00,651.00)(97.00,654.00)
\bezier17(185.00,581.00)(175.00,578.00)(168.00,581.00)
\bezier14(153.00,582.00)(143.00,582.00)(140.00,585.00)
\bezier14(126.00,588.00)(116.00,588.00)(113.00,591.00)
\bezier8(101.00,595.00)(96.00,597.00)(94.00,600.00)
\end{picture}
\end{center}
    \caption{The disks for the $3$-ended and $4$-ended Costa surfaces}
    \label{fig:costafund}
\end{figure}
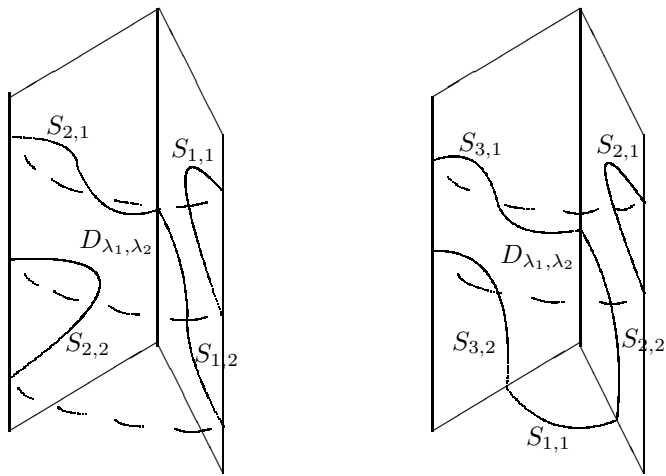

  The fundamental piece $D_{\lambda_1,\lambda_2}$ of this $4$-ended
  Costa example is shown in Figure~\ref{fig:costafund}.
  The boundary curves $S_{1,1},S_{2,1},S_{2,2},S_{3,1},S_{3,2}$ are 
  planar geodesics.  
  The planes $P_{2,1}$ and $P_{2,2}$ are parallel to the $x_1x_3$-plane, and
  $P_{3,1}$ and $P_{3,2}$ are parallel to the $x_2x_3$-plane.
  We set $P_{2,1}=\{x_2=c_1\}$,
        $P_{2,2}=\{x_2=c_2\}$,
        $P_{3,1}=\{x_1=c_3\}$,
  and
        $P_{3,2}=\{x_1=c_4\}$.  
  Then the functions $c_1$, $c_2$, $c_3$, and $c_4$ are continuous
  in $\lambda_1$ and $\lambda_2$, and 
  $f_1=c_2-c_1$ and $f_2=c_4-c_3$ are non-degenerate 
  \cite{Wohlgemuth}.  
  Thus the period problems are non-degenerate and by 
  Theorem~\ref{mainthm},
  we get the 4-ended Costa cousin.

It seems to be an 
interesting problem to construct 3-ended Costa surfaces.
Unlike the 4-ended Costa cousin, the fundamental piece of
the 3-ended Costa cousin contains a straight line in its boundary. Thus,
like the case of Enneper's surface, we must use the union of 
two fundamental
disks in our construction. Even though the period problem is nondegenerate
on the 
fundamental disk of a minimal 3-ended Costa surface, this does not imply 
the period problem is nondegenerate on the union of 
two fundamental disks.
However, existence of a 3-ended Costa surface has been verified numerically. 
(Clearly, this CMC-1 surface has less symmetries than the minimal 3-ended 
Costa surface.)

\section{Long time deformations of genus-0 CMC-$c$ surfaces}%
    \label{sec:long} 

The examples in the previous section and in Section~\ref{sec:fence}
are obtained through a small perturbation of $c$,
and in general we cannot explicitly determine the range of $c$ such that 
the corresponding \cmcc\ surfaces exist.
In this section, we find that we can determine the range under a certain
situation less general than before.  
We now assume that we have at most 
three smooth planar geodesics in the boundaries of the fundamental
pieces.
So we will use the simpler notation $\mu_j$ instead of the notation
$\mu_{j,i}$.

\subsection{Restricted situation}

\begin{figure}
\begin{center}
\unitlength=0.5pt
\begin{picture}(308.00,295.00)(147.00,410.00)
\put(296.00,670.00){\line(1,0){76.00}}
\put(296.00,670.00){\line(3,-5){7.00}}
\put(402.00,705.00){\vector(-3,-1){94.00}}
\put(332.00,439.00){\vector(-1,4){26.00}}
\put(323.00,621.00){\makebox(0,0)[cc]{$D$}}
\put(366.00,618.00){\makebox(0,0)[cc]{\footnotesize $S_1$}}
\put(338.00,548.00){\makebox(0,0)[cc]{\footnotesize $S_2$}}
\put(288.00,610.00){\makebox(0,0)[cc]{\footnotesize $S_3$}}
\put(412.00,700.00){\makebox(0,0)[cc]{$\frac{\pi}{2}$}}
\put(330.00,425.00){\makebox(0,0)[cc]{$\frac{\pi}{3}$}}
\bezier297(374.00,669.00)(301.00,538.00)(449.00,539.00)
\bezier297(151.00,540.00)(299.00,541.00)(226.00,671.00)
\bezier300(375.00,410.00)(300.00,540.00)(225.00,410.00)
\bezier25(364.00,665.00)(352.00,660.00)(340.00,662.00)
\bezier25(237.00,665.00)(244.00,659.00)(260.00,660.00)
\bezier36(286.00,657.00)(311.00,654.00)(322.00,657.00)
\bezier31(237.00,677.00)(253.00,682.00)(268.00,682.00)
\bezier39(288.00,684.00)(310.00,686.00)(327.00,683.00)
\bezier18(348.00,681.00)(361.00,679.00)(365.00,676.00)
\bezier33(149.00,534.00)(147.00,515.00)(150.00,501.00)
\bezier36(157.00,485.00)(163.00,466.00)(176.00,454.00)
\bezier31(190.00,436.00)(203.00,416.00)(211.00,414.00)
\bezier25(223.00,423.00)(221.00,435.00)(213.00,446.00)
\bezier26(158.00,533.00)(169.00,525.00)(176.00,514.00)
\bezier40(187.00,499.00)(202.00,483.00)(209.00,465.00)
\bezier29(450.00,533.00)(455.00,522.00)(449.00,506.00)
\bezier38(444.00,486.00)(439.00,465.00)(426.00,454.00)
\bezier37(413.00,433.00)(401.00,413.00)(387.00,411.00)
\bezier159(300.00,540.00)(408.00,461.00)(434.00,466.00)
\bezier72(300.00,540.00)(256.00,579.00)(243.00,573.00)
\bezier122(300.00,540.00)(213.00,484.00)(195.00,491.00)
\bezier66(300.00,540.00)(342.00,569.00)(357.00,567.00)
\bezier120(300.00,540.00)(294.00,647.00)(304.00,656.00)
\bezier69(300.00,540.00)(305.00,478.00)(298.00,475.00)
\bezier13(300.00,551.00)(307.00,551.00)(307.00,545.00)
\bezier20(346.00,585.00)(339.00,576.00)(336.00,567.00)
\bezier47(343.00,600.00)(325.00,582.00)(311.00,564.00)
\bezier63(350.00,632.00)(329.00,605.00)(304.00,589.00)
\bezier40(331.00,653.00)(318.00,645.00)(303.00,624.00)
\bezier11(299.00,665.00)(307.00,667.00)(306.00,670.00)
\end{picture}
\end{center}
        \caption{Assumption~\protect\ref{ass:genus0} ($m=2$, $n=3$)}
        \label{fig:assumption}
\end{figure}
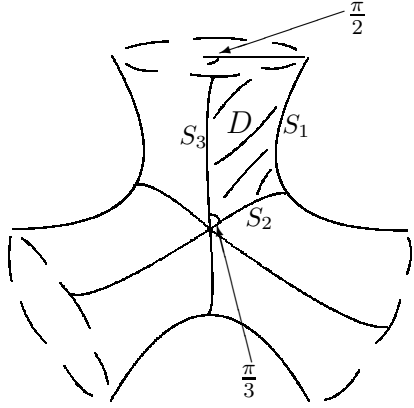

  Consider a complete minimal immersion $f_0\colon{}M\to \R^3$ of 
  a Riemann surface $M$.
  Throughout this section, we assume that $f_0$ is symmetric with
  respect to a disk $D\subset f(M)$ in the sense of
Definition~\ref{def:symmetric}
  and almost embedded in the sense of Definition~\ref{def:almost}.
  Moreover, we assume the following
  (cf.~Figure~\ref{fig:assumption})

  \begin{assum}\label{ass:genus0}
    The boundary $\partial D$ consists of three non-straight 
    planar geodesics $S_1$, $S_2$ and $S_3$ contained 
    in the planes $P_1$, $P_2$ and $P_3$, respectively, and
    \begin{itemize}
     \item   $S_1$ and $S_3$ are infinite rays, and $S_2$ is a
             curve of finite length;
     \item   $S_1$ and $S_2$ meet at $p_1$ with angle $\pi/n$,
             where $n\geq 3$ is an integer;
     \item   $S_2$ and $S_3$ meet at $p_2$ with angle $\pi/2$;
     \item   $S_1$ and $S_2$ bound a $1/(2m)$-piece of a catenoid end
             which is asymptotic to the standard catenoid,
             where $m\geq 2$ is an integer;
     \item   $1/2 < 1/m + 1/n$.
    \end{itemize}
   \end{assum}
  Here, the {\it standard catenoid\/} is the catenoid in $\R^3$ 
  obtained from the Gauss map $G_{\text{std}}=z$ and 
  the Hopf differential $Q_{\text{std}}=z^{-2}(dz)^2$.

  \begin{remark}\label{rem:alemb}
    Since $f_0$ is almost embedded, the planes $P_1$, $P_2$ and $P_3$ 
    are in general position (cf. Lemma~\ref{lem:threeplanes}).
    Hence we have:
    \begin{itemize}
      \item   Each pair of planes among  $P_1$, $P_2$ and $P_3$ 
              are not parallel.
      \item   Each pair of lines among $P_1\cap P_2$, $P_2\cap P_3$ and
              $P_3\cap P_1$ are not parallel.
    \end{itemize}
  \end{remark}        
  We have the following theorem.

  \begin{theorem}\label{thm:genus0gen}
    Let $f_0\colon{}M\to \R^3$ be a conformal minimal 
    immersion satisfying Assumption~\ref{ass:genus0},
    which is almost embedded, 
    and has Gauss map $G$ and Hopf differential $Q$.
        
    Then for each $c\neq 0$ satisfying
    \begin{equation}\label{eq:genus0cond}
       -\frac{1}{4}
          \left(
             \left(2-\frac{m}{2}+\frac{m}{n}\right)^2-1
          \right)<c<0
          \quad\text{or}\quad
      0<c<\frac{1}{4}
         \left(
           1-\left(\frac{m}{2}-\frac{m}{n}\right)^2
         \right),
    \end{equation}
    there exists an irreducible complete conformal \cmcc\ immersion
    $f_c\colon{}M\to \hypc$
    whose hyperbolic Gauss map and Hopf differential are $G$ and $Q$, 
    respectively.
  \end{theorem}

  We give a proof in Section~\ref{sub:proofgenus0}. 
  (About the interval
   \eqref{eq:genus0cond}, see Remark~\ref{rem:nec}.)

\begin{remark}\label{rem:ta}
Let $TA$ denote the total absolute curvature
$\int_M(-K)dA$,  
where $dA$ is the area
  element with respect to the induced metric $ds^2$.
When $c$ varies over the range \eqref{eq:genus0cond},
$TA$ is given by
  \[
     TA = 2\pi\left[N(\sqrt{1-4c}-1)+2N-2\right].
  \]
  where $N$ denotes the number of the ends.
  This is verified as follows:
  The total absolute curvature is the area with
  respect to the pseudo-metric $d\sigma_f^2$ in \eqref{eq:Def}.
  Thus, by Gauss-Bonnet formula, we have
  ${TA}/{2\pi}=\chi(\overline M)+\sum_{p\in\overline M}\ord_p {d\sigma_f^2}$.
  Moreover, by \eqref{eq:metHop}, $\ord_p d\sigma_f^2=\ord_p Q$ on $M$
  because $ds^2$ is non-degenerate.
  Thus, 
  \[
     \frac{TA}{2\pi}=\chi(\overline M) + 
                        \sum_{p \in \overline M\setminus M} \ord_p d\sigma_f^2
                        + \sum_{p \in M}\ord_p Q.
  \]
  Since $\overline M$ is of genus $0$ and $Q$ has poles of order $2$ at
  the ends, the sum of $\ord_p Q$   over $p\in M$ is $2N-4$.
  On the other hand, by the upcoming equation \eqref{eq:genus0expanda}, 
  $\ord_p d\sigma_f^2$ is $\lambda-1$, where $\lambda=\sqrt{1-4c}$.
  Hence, we have the formula. We
  remark that the total absolute curvature of the corresponding minimal
  surface in $\R^3$ is $4\pi(N-1)$.
  Thus $TA>4\pi(N-1)$ (resp.~$TA<4\pi(N-1)$) if $c<0$ (resp.~$c>0$).
\end{remark}

\subsection{Examples}
  The following two examples
  have already  been constructed in \cite{UY4}.  
  However the method in this section is more explicit and suggests 
  an algorithm to draw the surfaces by numerical calculations.

  \begin{example}[Genus-$0$ Jorge-Meeks surface]
    \label{ex:genus0jm}
    Let $f_0\colon{}M\to\R^3$ be the Jorge-Meeks $n$-oid, 
    with each end asymptotic to a standard catenoid.
    Then the fundamental region $D$ of $f_0$ satisfies 
    Assumption~\ref{ass:genus0}, where $m=2$ and $n$ is the number of ends.
    Then by Theorem~\ref{thm:genus0gen}, 
    for each $c$ satisfying
    $-(n+1)/{n^2}<c<0$ or $0<c<(n-1)/{n^2}$,
    there exists a conformal \cmcc\ immersion $f_c\colon{}M\to \hypc$
    whose hyperbolic Gauss map and Hopf differential coincide 
    with those of $f_0$.
    The total absolute curvature $TA$ of $f_c$ varies over the range
    $\left(\mathstrut 4\pi(n-2),4\pi(n-1)\right)\cup
     \left(\mathstrut 4\pi(n-1),4\pi n)\right)$ (cf.~Remark~\ref{rem:ta}).

\begin{figure}
\vspace{2in}
    \caption{Profile curves for the genus-0 Jorge-Meeks cousin 
      ($c<0$ on the left, $c>0$ on the right).}   
    \label{fig:jm2}    
\end{figure}    

  \end{example}

  \begin{example}[Genus-$0$ surfaces with Platonic symmetry]
  \label{ex:gen0platonic}
    There are minimal surfaces with symmetries of the other Platonoids.
    For such minimal surfaces, we can  apply Theorem~\ref{thm:genus0gen}.
    In these cases, $n$ is the number of 
    edges of the Platonic solid bounding each face, and
    $m$ is the number of edges with a common vertex.
    Table~\ref{tab:platonoid}  shows the range of $c$ for which we know the 
    minimal surfaces can be deformed to \cmcc{} surfaces.
\begin{table}[htb]
    \begin{center}
    \begin{tabular}{l|c|c|c|c|c}
                 & Number   &      &     &            &  \\
      Symmetry   & of ends  &  $m$ & $n$ &  Range of $c$ &Range of $TA$\\
\hline
      Tetrahedra & $4$      &  $3$ & $3$ &  
                   $\left(-\frac{\mathstrut{}5}{\mathstrut{}16},0\right)\cup
                    \left(0,\frac{\mathstrut{}3}{\mathstrut{}16}\right)$&
                   $(8\pi,12\pi)\cup(12\pi,16\pi)$ \\ \hline
      Hexahedra  & $8$      &  $3$ & $4$ &  
                   $\left(-\frac{\mathstrut{}9}{\mathstrut{}64},0\right)\cup
                    \left(0,\frac{\mathstrut{}7}{\mathstrut{}64}\right)$&
                   $(24\pi,28\pi)\cup(28\pi,32\pi)$ \\ \hline
      Octahedra  & $6$      &  $4$ & $3$ &  
                   $\left(-\frac{\mathstrut{}7}{\mathstrut{}36},0\right)\cup
                    \left(0,\frac{\mathstrut{}5}{\mathstrut{}36}\right)$&
                   $(16\pi,20\pi)\cup(20\pi,24\pi)$ \\ \hline
      Dodecahedra& $20$     &  $3$ & $5$ &  
                   $\left(-\frac{\mathstrut{}21}{\mathstrut{}400},0\right)\cup
                    \left(0,\frac{\mathstrut{}19}{\mathstrut{}400}\right)$&
                   $(72\pi,76\pi)\cup(76\pi,80\pi)$ \\ \hline
      Icosahedra & $12$     &  $5$ & $3$ &  
                   $\left(-\frac{\mathstrut{}13}{\mathstrut{}144},0\right)\cup
                    \left(0,\frac{\mathstrut{}11}{\mathstrut{}144}\right)$& 
                   $(40\pi,44\pi)\cup(44\pi,48\pi)$
    \end{tabular}
    \end{center}
\caption{}\label{tab:platonoid}
\end{table}
  \end{example}
\begin{remark}
    The corresponding minimal surfaces of \cmcc\ surfaces in 
    Example~\ref{ex:genus0jm} and \ref{ex:gen0platonic} are non-embedded.
    Although $f_c$ has self-intersections for $c<0$,
    it seems to be embedded for some positive $c$ (the total absolute
    curvature is less than that of the corresponding minimal surface),
    by our numerical calculations  (see Figure~\ref{fig:jm2}).
\end{remark}
\subsection{Proof of Theorem~\protect\ref{thm:genus0gen}}%
          \label{sub:proofgenus0}
  Let $\mu_j$ be the reflection across the plane $P_j$ containing 
  the geodesic $S_j$.
  In the same way as in the previous section, 
  we can induce the reflections $\tilde\mu_j$ on
  the universal cover $\widetilde M$.
  Let $l$ be a loop on $M$ surrounding the end 
  that intersects with $D$, 
  and let $\tau\in\pi_1(M)$ be 
  the deck transformation on $\widetilde M$ induced by $l$.
  Then the reflections satisfy the following relations:
  \begin{equation}\label{eq:genus0refrel}
      \tilde\mu_j\circ\tilde\mu_j=\id\,(j=1,2,3),\quad
                (\tilde\mu_1\circ\tilde\mu_2)^n=\id,\quad
                (\tilde\mu_2\circ\tilde\mu_3)^2=\id,\quad
                (\tilde\mu_1\circ\tilde\mu_3)^m=\tau.
  \end{equation}
  The Gauss map $G$ and the Hopf differential $Q$ of $f$ can be lifted
  to $\widetilde M$, where we will continue to denote them 
  by $G$ and $Q$, respectively.
  Then, by symmetry, 
  there exist matrices $\sigma(\mu_j)\in SU(2)$ ($j=1,2,3$) such that
  \begin{equation}\label{eq:genus0gaussref}
        \overline{G\circ\tilde\mu_j}=\sigma(\mu_j)^{-1}* G,\qquad
        \overline{Q\circ\tilde\mu_j}=Q
        \qquad (j=1,2,3).   
  \end{equation}
  \begin{lemma}\label{lem:genus0sigma}
    By a suitable choice of the coordinate system of $\R^3$, we can choose
    \begin{equation}\label{eq:genus0sigma}
        \sigma(\mu_1)= \id,\quad
        \sigma(\mu_2)=
             \begin{pmatrix}
                e^{{\pi i}/{n}}  &   0   \\
                0   &   e^{-{\pi i}/{n}}
             \end{pmatrix},
                \qquad
        \sigma(\mu_3)=i
             \begin{pmatrix}
                \alpha_0 e^{{\pi i}/{n}}   &   \beta_0 \\
                \beta_0   &   -\alpha_0 e^{-{\pi i}/{n}}
             \end{pmatrix}, 
    \end{equation}
    where 
    $\alpha_0=\cos({\pi}/{m})/\sin({\pi}/{n}) \in \R$
    and $\alpha_0^2+\beta_0^2=1$, $\beta_0 \in \R$.
  \end{lemma}
  \begin{pf}
    Take a coordinate system of $\R^3$ 
    such that the plane $P_1$ is the $x_1x_3$-plane and 
    $P_1\cap P_2$ is the $x_3$-axis.
    Since $P_1$ and $P_2$ form an angle $\pi/n$, 
    $\sigma(\mu_1)$ and $\sigma(\mu_2)$ satisfies 
     \eqref{eq:genus0sigma}.
    By $\tilde\mu_3\circ\tilde\mu_3=\id$ and Lemma~\ref{lem:id},
    we have ${\sigma(\mu_3)}\overline{\sigma(\mu_3)}=\id$.
    This shows that 
    \begin{equation}\label{eq:genus0sigma3}
         \sigma(\mu_3)=
            \begin{pmatrix}
                a   &   i\beta_0 \\
                i\beta_0    &   \bar a
            \end{pmatrix}\qquad
           (a\bar a+\beta_0^2=1, \, \beta_0 \in \R) \; .
    \end{equation}

    By the relation $(\tilde\mu_2\circ\tilde\mu_3)^2=\id$,
    we have $({\sigma(\mu_2)}\overline{\sigma(\mu_3)})^2=\pm\id$.
    If $({\sigma(\mu_2)}\overline{\sigma(\mu_3)})^2=\id$,
    by \eqref{eq:genus0sigma3} and the form of $\sigma(\mu_2)$,
    we have $a=\pm e^{-{\pi i}/{n}}$ and $\beta_0=0$.
    This shows that  $\sigma(\mu_3)=\sigma(\mu_2)$,
    and that $P_2$ and $P_3$ are parallel.
    This is impossible by Remark~\ref{rem:alemb}.
    Hence
    \begin{equation}\label{eq:genus0sigma23}
          \left({\sigma(\mu_2)}\overline{\sigma(\mu_3)}\right)^2=-\id
           \; \; ,
    \end{equation}
    and then 
    \begin{equation}\label{eq:genus0sigma3a}
           \sigma(\mu_3)=i
            \begin{pmatrix}
               \alpha_0 e^{{\pi i}/{n}}    &   \beta_0 \\
               \beta_0    &   -\alpha_0 e^{-{\pi i}/{n}}
             \end{pmatrix}\qquad
                (\alpha_0^2+\beta_0^2=1) \; \; ,
    \end{equation}
    where $\alpha_0$ and $\beta_0$ are real numbers.
    Since the rays $S_1$ and $S_3$ span a $1/(2m)$-piece of a 
    catenoid end, 
    the angle between the planes $P_1$ and $P_3$ is $\pi/m$.
    Hence the eigenvalues of $\sigma(\mu_3)$ are 
    $e^{\pm {\pi i}/{m}}$.
    Thus, we have 
        $\trace\sigma(\mu_3)=2\cos(\pi/{m})
                =2\alpha_0\sin({\pi}/{n})$.
  \end{pf}

  The metric $ds_G^2$ defined in 
  Lemma~\ref{lem:Dempty} is positive definite on $M$.
  Hence the set $\D_{M}^{(c)}(G,Q)$ is not empty.

  \begin{proposition}\label{prop:genus0step2}
    For each $c$ satisfying \eqref{eq:genus0cond},
    there exists a one-parameter family of
    null holomorphic immersion
    of $\widetilde{M}$
    into $SL(2,\C)$ with the following properties.
    \begin{enumerate}
     \item   \label{item:1}
             $f_{c}=(1/c)F_{c}{F_{c}^*}\in
             \D_{M}^{(c)}(G,Q)$ 
             for each $c$.
     \item   \label{item:2}
             $F_{c}$ is smooth in $c$.
     \item   \label{item:genus00}
             $\lim_{c\to 0}F_{c}=\id$.
     \item   $\lim_{c\to 0}\hat\rho^{}_{F_{c}}(\tilde\mu_{j})=
              \sigma(\mu_{j})$\,\,$(j=1,2,3)$.
     \item   \label{item:genus01}
             $\hat\rho^{}_{F_{c}}(\tilde\mu_{1})=\id$.
     \item   \label{item:genus02}
             There exists real-valued continuous functions $\alpha$ 
             and $\beta$ in $c$ such that 
             \[
                 \hat\rho_{F_c}(\tilde\mu_2)=
                    \begin{pmatrix}
                        e^{{\pi i}/{n}} & 0  \\
                        0 & e^{-{\pi i}/{n}}
                    \end{pmatrix},
                 \quad
                 \hat\rho_{F_c}(\tilde\mu_3)
                 = i\begin{pmatrix}
                        \alpha e^{{\pi i}/{n}} &   \beta\\
                        \beta  &   -\alpha e^{-{\pi i}/{n}}
                    \end{pmatrix}\in SU(2),
              \]
    \end{enumerate}
    $\alpha^2+\beta^2=1$, and $\lim_{c\to 0}\alpha(c)=\alpha_0$.
  \end{proposition}

  Once the above proposition is proven, $f_c\in I_M^{(c)}(G,Q)$
  by Proposition~\ref{prop:su2cond}. 
  Thus $f_c$ is (after an appropriate `dilation') the desired
  \cmcone{} immersion as in Theorem~\ref{thm:genus0gen}.

  \begin{pf*}{\it Proof of Proposition~\ref{prop:genus0step2}}
    In the same way as in Step I of the proof of
    Proposition~\ref{prop:JMStep2}, 
    we can choose $\{\check F_c\}$ which satisfies \ref{item:genus00} and 
    \ref{item:genus01}.

    The eigenvalues of
    $\hat\rho_{\check F_c}(\tilde\mu_2)$ are 
    $\{e^{{\pi i}/{n}},e^{{-\pi i}/{n}}\}$.
    Indeed, since
    $(\overline{\sigma(\mu_1)}{\sigma(\mu_2)})^n=-\id$ and 
    $\hat\rho_{\check F_c}(\tilde\mu_1)=\id$,
    we have $\hat\rho_{\check F_c}(\tilde\mu_2)^n=-\id$.
    Then the eigenvalues of $\hat\rho_{\check F_c}(\tilde\mu_2)$ are
    $n$-th roots of $-1$.
    Then, by the continuity of $\hat\rho_{\check F_c}(\tilde\mu_2)$,
    they must be constant.
    Here, $\lim_{c\to 0}\hat\rho_{\check F_c}(\tilde\mu_2)=\sigma(\mu_2)$,
    so the eigenvalues of $\hat\rho_{\check F_c}(\tilde\mu_2)$ coincide
    with those of $\sigma(\mu_2)$.
    In particular, the eigenvalues 
    of $\hat\rho_{\check F_c}(\tilde\mu_2)$ are not real.
    Then, by the same method as in Step II of the proof
    of Proposition~\ref{prop:JMStep2},
    we can diagonalize it with a real matrix $u=u(c)$:
    \[
        u(c)^{-1}\hat\rho_{\check F_c}^{}(\tilde\mu_2)u(c)
            =\begin{pmatrix}
                e^{{\pi i}/{n}}  & 0 \\
                0 & e^{{-\pi i}/{n}}
             \end{pmatrix}.
    \]
    Setting  $\hat F_c=\check F_c u(c)$, we obtain a family $\{\hat F_c\}$
    satisfying \ref{item:genus00}, \ref{item:genus01} and 
    \ref{item:genus02}.
    Moreover, by using $(\tilde\mu_2\circ\tilde\mu_3)^2=\id$ and 
    \eqref{eq:genus0sigma23},
    $\left({\hat\rho_{\hat F_c}(\tilde{\mu}_2)}
     \overline{\hat\rho_{\hat F_c}(\tilde{\mu}_3)}\right)^2
            =-\id$
    holds and then $\hat\rho_{\hat F_c}(\tilde\mu_3)$ can be written
    in the form
    \begin{equation}\label{eq:genus0rho3}
        \hat\rho_{\hat F_c}(\tilde\mu_3)
          =i\begin{pmatrix}
                e^{{\pi i}/{n}}\alpha   &   \beta_1    \\
                \beta_2    & -e^{-{\pi i}/{n}}\alpha
            \end{pmatrix}
            \qquad
            \left(\alpha,\beta_1,\beta_2\in\R,\quad
            \alpha^2+\beta_1\beta_2=1\right). 
    \end{equation}

    When $c$ is sufficiently small, 
    we can prove the proposition
    by applying the
    same argument as in Step III of
    the proof of Proposition~\ref{prop:JMStep2},
    because of $\beta_1/\beta_2>0$.
    So it is sufficient to show that
    $\beta_1/\beta_2>0$ for any $c$ satisfying \eqref{eq:genus0cond}.
    Thus the proof of the proposition reduces to the following lemma.
  \end{pf*}
 
  \begin{lemma}\label{lem:genus0step3}
    If $c$ satisfies \eqref{eq:genus0cond},
    then $\beta_1\beta_2>0$.
  \end{lemma}
  \begin{pf}
    Since $\sigma(\mu_3)$ has plus-minus ambiguity, we have two
    choices for $\sigma(\mu_3)$.  We choose $\sigma(\mu_3)$ so that 
    the eigenvalues of $\sigma(\mu_3)$ are $e^{\pm {\pi i}/{m}}$,
    and then 
    \[
         \left({\sigma(\mu_1)}\overline{\sigma(\mu_3)}\right)^m
            =\sigma(\mu_3)^m=-\id \; \; .
    \]
    By the last relation in \eqref{eq:genus0refrel}, we have
    \begin{equation}\label{eq:genus0deck1}
         \rho_{F_c}(\tau)=-\left({\hat\rho_{F_c}({\tilde\mu}_1)}
                           \overline{\hat\rho_{F_c}({\tilde\mu}_3)} \right)^m
                         =-\overline{\hat\rho_{F_c}({\tilde\mu}_3)}^m,
    \end{equation}
    where $\tau$ is the deck transformation 
    induced from the loop $l$ surrounding the catenoid end.

    The eigenvalues of $\hat\rho_{F_c}(\tilde{\mu}_3)$ are 
    $\xi_\pm=s\pm i\sqrt{1-s^2}$, where
    $s=\alpha\sin(\pi/n)$.
    Thus, if
    \begin{equation}\label{eq:alpharange}
        |\alpha|<\frac{1}{\sin({\pi}/{n})} \; \; ,
    \end{equation}
    $\xi_\pm$ are complex numbers of absolute value $1$, 
    and then, there exists $\theta\in(0,\pi)$ such that 
    $\xi_\pm=e^{\pm i\theta}$ and 
    $\theta=\cos^{-1}\left(\alpha\sin({\pi}/{n})\right)$.
    Since $\alpha(0)=\alpha_0$ satisfies \eqref{eq:alpharange},
    there exists an interval $I_0\subset \R$ which contains the origin
    and \eqref{eq:alpharange} holds for each $c\in I_0$.  
    We now restrict $c\in I_0$ in order to make $|\xi_k|=1$.
    Later, we restrict the range of $c$ again to make
    $|\alpha(c)|<1$.
    Hence, for each $c\in I_0$, it holds that
    \begin{equation}\label{eq:tracerho}
        \trace \rho_{F_c}(\tau)=-
        \overline{\trace\left(\hat\rho_{F_c}(\tilde\mu_3)\right)}^m
                 =-(e^{-im\theta}+e^{im\theta})=-2\cos m\theta.
    \end{equation}
        
    Let $\overline M$ be the compactification of $M$ and $p\in \overline M$
    a point corresponding to an end.
    Since the end is asymptotic to a standard catenoid, 
    the Weierstrass data are also asymptotic to those of a catenoid.  
    Hence there exists a coordinate $z$ of $\overline M$ such that
    $z(p)=0$ and $(G,Q)$ are expanded as
    \begin{equation}\label{eq:genus0expand}
         G=z+\dots,\qquad Q=\left(\frac{1}{z^2}+\dots\right)(dz)^2.
    \end{equation}

    Let $g_c$ be the secondary Gauss map of $F_c$.
    Then, by (\ref{eq:pull}), we have
    \begin{equation}\label{eq:genus0gausschange}
        g\circ \tau=\rho_{F_c}(\tau)* g
                   =\frac{a_{11}g+a_{12}}{a_{21}g+a_{22}}\qquad
        \text{where  }\rho_{F_c}(\tau)=(a_{kj}).
    \end{equation}
    So, by \eqref{eq:S} and the fact that the Schwarzian derivative
    is invariant under M\"obius transformations, 
    there exists a matrix $b\in SL(2,\C)$ such that 
    \begin{equation}\label{eq:genus0expanda}
            b* g=
                z^{\lambda}(g_0+g_1 z+g_2 z^2+\dots),\qquad
                    g_0\neq 0,
    \end{equation}
    where
    $\lambda= \sqrt{1-4c}$.
    Hence, by \eqref{eq:genus0gausschange},
    we have 
    \[
         b^{-1}\rho_{F_c}(\tau)b
                    =\pm\begin{pmatrix}
                        e^{\pi\lambda i} & 0 \\
                        0  & e^{-\pi\lambda i}
                     \end{pmatrix}\;\; ;
    \]
    that is, 
    the eigenvalues of $\rho_{F_c}(\tau)$ 
    are $\{\pm e^{\pi\lambda i},\pm e^{-\pi\lambda i}\}$.
    So, for each $c<{1}/{4}$, 
    \begin{equation}\label{eq:tracerho1}
         \trace\rho_{F_c}(\tau)=\pm 2\cos\pi\lambda
    \end{equation}
    holds.
    Combining \eqref{eq:tracerho} and \eqref{eq:tracerho1},
    we have
    $\pm\cos\pi\lambda=-\cos m \theta$.
    Letting $c\to 0$, $\theta$ tends to 
    $\theta(0)=\pi/m$ and $\lambda\to 1$.
    Thus the sign on the left-hand side is ``minus":
    \begin{equation}\label{eq:lambdatheta}
         \cos\pi\lambda=\cos m \theta.
    \end{equation}
        
    Now we assume $c$ satisfies
    \eqref{eq:genus0cond}, that is, $c \in I$, where $I\setminus\{0\}$ is the
    open interval 
    \begin{equation}\label{eq:cint}
         I:=\left(-\frac{1}{4}\left\{
            \left(2-\frac{m}{2}+\frac{m}{n}\right)^2-1\right\},
            \frac{1}{4}
            \left\{1-\left(\frac{m}{2}-\frac{m}{n}\right)^2\right\}\right)
                    \; \; . 
    \end{equation}
    So, if $c\in I$, we have 
    \[
        \pi m\left(({1}/{2})-({1}/{n})\right)  <m\theta
            <2\pi-\pi m\left(({1}/{2})-({1}/{n})\right),
    \]
    because of \eqref{eq:lambdatheta}.

    Here, $\theta=\theta(c)$ is a continuous function, 
    and $m\theta(0)=\pi$, and 
    $m\theta(I)$ is an interval containing $\pi$.
    Thus, we have
    \[
         \theta(I)\subset \left(
         \pi\left(\frac{1}{2}-\frac{1}{n}\right),
         \pi\left(\frac{2}{m}-\frac{1}{2}+\frac{1}{n}\right)\right)
         \subset (0,\pi).
    \]
    On $(0,\pi)$, $\cos\theta$ is a decreasing function.
    Thus 
    \[
         \sin\frac{\pi}{n} =
         \cos\left(\pi\left(\frac{1}{2}-\frac{1}{n}\right)\right)
                >\cos\theta(c)>
         \cos\left(\pi\left(\frac{2}{m}-\frac{1}{2}+\frac{1}{n}\right)\right)
         \geq -\sin\frac{\pi}{n} \; \; .
    \]
    This shows
    \[
        -\sin({\pi}/{n})<\cos\theta(c)=\alpha\sin({\pi}/{n})
          <\sin({\pi}/{n}).
    \]
    Therefore, $|\alpha(c)|<1$ if $c\in I$.
    For such a value of $c$,
    $\beta_1\beta_2=1-\alpha^2>0$.
  \end{pf}

  \begin{remark}\label{rem:nec}
    The interval (\ref{eq:genus0cond}) in Theorem~\ref{thm:genus0gen}
    gives only a sufficient
    condition for the existence of deformation.
    However, in the case of the genus $0$ Jorge-Meeks surfaces, 
    one can determine a necessary and sufficient condition
    on $c$ so that there exists a corresponding \cmcc{} surface.
    Indeed, for the Jorge-Meeks $n$-oid $f_0\colon{}M\to\R^3$ 
    as in Example~\ref{ex:genus0jm}, 
    there exists a corresponding \cmcc{} immersion
    into $\hypc$ {\em if and only if\/} $c\in I_{0,+}$, $I_{0,-}$
    or $I_k$ ($k=1,2,\dots$), where 
    \begin{equation}\label{eq:jmnecsuf}
       \begin{aligned}
           I_{0,+} &= \left(0,\frac{n-1}{n^2}\right),\qquad
           I_{0,-} = \left(-\frac{n+1}{n^2},0\right)\\
           I_{k\hphantom{,-}}   &=
                      \left(
                      -\left(k+\frac{1}{n}\right)\left(k+1+\frac{1}{n}\right),
                      \left(k-\frac{1}{n}\right)\left(k+1-\frac{1}{n}\right)
                      \right)
        \quad (k=1,2,\dots).
       \end{aligned}
    \end{equation}   
 Since $m=2$ in this case, 
    the condition (\ref{eq:lambdatheta}) can be rewritten as
    $\cos\pi\lambda=\cos 2\theta=2\cos^2\theta-1$,
    and recalling that $\cos\theta=\alpha\sin(\pi/n)$
    (see the first part in the proof of Lemma~\ref{lem:genus0step3}),
    we have
    $\cos\pi\lambda=2\alpha^2\sin^2{({\pi}/{n})}-1$,
    and then,
    $\alpha^2=(\cos\pi\lambda+1)/(2\sin^2{({\pi}/{n})})$
    holds.
    Here, a corresponding \cmcc{} immersion exists if and only if 
    the conclusion of Lemma~\ref{lem:genus0step3}, 
    i.e., $\beta_1\beta_2=1-\alpha^2>0$ holds.
    (Otherwise, one cannot find $f_c$ with $SU(2)$-condition.)
    Using the above relation,
    this necessary and sufficient condition is rewritten as
    $\cos\pi\lambda+1< 2\sin^2({\pi}/{n})$,
     where 
    $\lambda=\sqrt{1-4c}$.
    This inequality holds if and only if \eqref{eq:jmnecsuf} holds.

    By the same argument as Remark~\ref{rem:ta}, the total absolute
    curvature for $f_c$ ($c\in I_k$) varies over the interval 
    $\left(\mathstrut 4\pi\{n(k+1)-2\},4\pi n(k+1)\right)$.
    Numerical investigations suggest that 
the surface is embedded for $c\in I_{0,+}$ sufficiently near $(n-1)/n^2$. 
  \end{remark}

\section{Periodic \cmcc\ surfaces}\label{sec:fence}

The construction in Section~\ref{sec:ftc}
depends on the properties of minimal surfaces in $\R^3$.
However, our method can sometimes be applied even 
if the corresponding minimal surfaces do not exist.  
We believe that the following example is of interest, 
because our construction method here is independent of 
the corresponding minimal surfaces.  

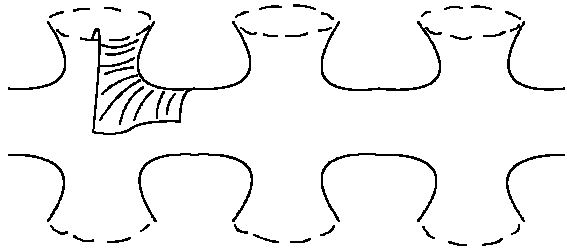
\begin{figure}
\begin{center}
\unitlength=0.5pt
\begin{picture}(420.00,189.00)(90.00,621.00)
\bezier130(90.00,740.00)(156.00,736.00)(120.00,790.00)
\bezier131(200.00,790.00)(165.00,734.00)(230.00,740.00)
\bezier124(90.00,690.00)(154.00,690.00)(120.00,640.00)
\bezier125(201.00,640.00)(167.00,693.00)(230.00,690.00)
\bezier130(230.00,740.00)(296.00,736.00)(260.00,790.00)
\bezier135(370.00,740.00)(303.00,732.00)(340.00,790.00)
\bezier134(230.00,690.00)(299.00,693.00)(260.00,640.00)
\bezier132(370.00,690.00)(303.00,694.00)(340.00,640.00)
\bezier129(510.00,740.00)(446.00,734.00)(480.00,790.00)
\bezier126(510.00,690.00)(444.00,689.00)(480.00,640.00)
\bezier131(370.00,690.00)(438.00,691.00)(400.00,640.00)
\bezier146(370.00,740.00)(444.00,733.00)(400.00,790.00)
\bezier10(127.00,784.00)(132.00,779.00)(135.00,780.00)
\bezier12(145.00,778.00)(153.00,776.00)(157.00,777.00)
\bezier15(167.00,777.00)(177.00,776.00)(181.00,779.00)
\bezier10(187.00,781.00)(194.00,784.00)(196.00,787.00)
\bezier12(197.00,793.00)(197.00,797.00)(189.00,798.00)
\bezier11(174.00,799.00)(170.00,801.00)(163.00,801.00)
\bezier11(151.00,800.00)(144.00,801.00)(140.00,799.00)
\bezier10(131.00,798.00)(126.00,797.00)(121.00,794.00)
\bezier13(266.00,786.00)(273.00,780.00)(277.00,782.00)
\bezier10(283.00,781.00)(289.00,779.00)(293.00,780.00)
\bezier16(300.00,780.00)(311.00,779.00)(316.00,782.00)
\bezier9(325.00,784.00)(331.00,784.00)(334.00,786.00)
\bezier9(335.00,793.00)(331.00,797.00)(327.00,797.00)
\bezier16(322.00,800.00)(316.00,803.00)(306.00,803.00)
\bezier16(296.00,802.00)(287.00,803.00)(280.00,800.00)
\bezier10(271.00,798.00)(267.00,798.00)(262.00,794.00)
\bezier13(406.00,786.00)(411.00,780.00)(417.00,782.00)
\bezier15(425.00,780.00)(434.00,776.00)(440.00,779.00)
\bezier16(450.00,780.00)(460.00,779.00)(465.00,783.00)
\bezier8(470.00,784.00)(476.00,786.00)(476.00,788.00)
\bezier12(477.00,793.00)(470.00,797.00)(466.00,798.00)
\bezier14(455.00,800.00)(448.00,803.00)(441.00,802.00)
\bezier13(433.00,801.00)(425.00,803.00)(420.00,800.00)
\bezier8(409.00,798.00)(405.00,799.00)(403.00,795.00)
\bezier15(122.00,636.00)(129.00,628.00)(134.00,631.00)
\bezier15(140.00,628.00)(147.00,623.00)(154.00,625.00)
\bezier17(165.00,624.00)(177.00,624.00)(182.00,627.00)
\bezier7(190.00,632.00)(195.00,635.00)(196.00,637.00)
\bezier13(261.00,637.00)(263.00,631.00)(270.00,631.00)
\bezier19(276.00,629.00)(288.00,624.00)(294.00,624.00)
\bezier14(304.00,624.00)(312.00,623.00)(318.00,626.00)
\bezier14(326.00,629.00)(335.00,631.00)(336.00,636.00)
\bezier12(404.00,635.00)(407.00,628.00)(412.00,627.00)
\bezier13(423.00,624.00)(432.00,621.00)(436.00,622.00)
\bezier15(449.00,622.00)(461.00,623.00)(463.00,626.00)
\bezier11(471.00,630.00)(477.00,632.00)(476.00,637.00)
\bezier34(230.00,740.00)(219.00,738.00)(220.00,715.00)
\bezier34(220.00,715.00)(200.00,716.00)(186.00,712.00)
\bezier35(186.00,712.00)(176.00,703.00)(154.00,707.00)
\bezier137(154.00,707.00)(164.00,810.00)(153.00,777.00)
\bezier20(207.00,737.00)(202.00,727.00)(202.00,718.00)
\bezier35(192.00,741.00)(179.00,734.00)(174.00,713.00)
\bezier30(187.00,750.00)(170.00,746.00)(161.00,736.00)
\bezier26(185.00,763.00)(172.00,756.00)(160.00,757.00)
\bezier30(188.00,777.00)(173.00,768.00)(161.00,773.00)
\bezier17(216.00,736.00)(210.00,730.00)(210.00,721.00)
\bezier26(199.00,738.00)(190.00,731.00)(185.00,716.00)
\bezier43(189.00,746.00)(173.00,738.00)(159.00,716.00)
\bezier21(184.00,756.00)(173.00,755.00)(163.00,751.00)
\bezier22(185.00,771.00)(175.00,765.00)(164.00,766.00)
\end{picture}
\end{center}
    \caption{The catenoid fence and its fundamental disk}%
    \label{fig:catfence}
\end{figure}

In this section, we shall construct singly-periodic \cmcc\ surfaces
in $\hypc$.
The corresponding minimal surfaces are well-known
as the ``catenoid fence" and ``Jorge-Meeks $n$-oid fence"
(cf.\ \cite{Karcher}, \cite{Rossman}, Figure~\ref{fig:catfence}).
Although, in the following argument, we borrow geometric intuition
from the corresponding minimal surfaces,
our construction will not depend on them.  
Thus, by letting $c\to 0$, we have another way to construct the 
minimal catenoid (resp. $n$-oid) fence.

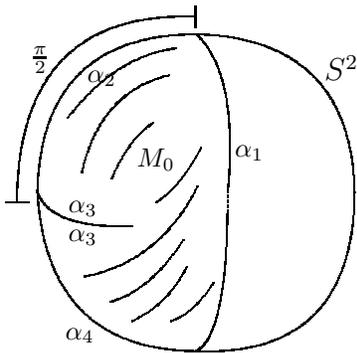
\begin{figure}
\begin{center}
\unitlength=0.6pt
\begin{picture}(220.00,222.00)(130.00,586.00)
\put(130.00,684.00){\line(1,0){15.00}}
\put(250.00,808.00){\line(0,-1){14.00}}
\put(225.00,711.00){\makebox(0,0)[cc]{\footnotesize $M_{0}$}}
\put(342.00,769.00){\makebox(0,0)[cc]{$S^2$}}
\put(152.00,774.00){\makebox(0,0)[cc]{\footnotesize $\frac{\pi}{2}$}}
\put(284.00,715.00){\makebox(0,0)[cc]{\footnotesize $\alpha_1$}}
\put(191.00,762.00){\makebox(0,0)[cc]{\footnotesize $\alpha_2$}}
\put(179.00,680.00){\makebox(0,0)[cc]{\footnotesize $\alpha_3$}}
\put(179.00,662.00){\makebox(0,0)[cc]{\footnotesize $\alpha_3$}}
\put(177.00,600.00){\makebox(0,0)[cc]{\footnotesize $\alpha_4$}}
\bezier197(250.00,790.00)(152.00,789.00)(150.00,690.00)
\bezier195(250.00,790.00)(348.00,787.00)(350.00,690.00)
\bezier203(350.00,690.00)(349.00,586.00)(250.00,590.00)
\bezier193(250.00,590.00)(156.00,591.00)(150.00,690.00)
\bezier111(250.00,790.00)(277.00,765.00)(270.00,690.00)
\bezier112(270.00,690.00)(270.00,600.00)(250.00,590.00)
\bezier75(150.00,690.00)(159.00,667.00)(210.00,669.00)
\bezier233(137.00,684.00)(138.00,805.00)(250.00,801.00)
\bezier96(233.00,764.00)(190.00,754.00)(178.00,703.00)
\bezier45(223.00,734.00)(207.00,723.00)(197.00,699.00)
\bezier101(251.00,694.00)(228.00,647.00)(180.00,637.00)
\bezier50(244.00,643.00)(233.00,620.00)(210.00,609.00)
\bezier45(253.00,718.00)(239.00,692.00)(225.00,684.00)
\bezier63(242.00,660.00)(221.00,628.00)(196.00,621.00)
\bezier37(261.00,633.00)(248.00,613.00)(234.00,609.00)
\bezier87(237.00,781.00)(199.00,776.00)(169.00,737.00)
\end{picture}
\end{center}
        \caption{The image of $M_0$ by $G_0$}%
        \label{fig:gaussim}
\end{figure}

We denote the shaded region in Figure~\ref{fig:catfence} by $M_0$.
It is the fundamental region of the catenoid fence,
that is, the whole surface is obtained by reflections 
of $M_0$ about its boundary curves.
In the \cmcc{} case, 
we define the underlying Riemann surface and the data $g$ and $Q$
by modifying $M_0$.
Let $G_0$ and $Q_0$ be the Gauss map and the Hopf differential
of the minimal surface $M_0$.
We consider $G_0$ as a map to the unit sphere $S^2$.
So the image of $G_0$ in $S^2$ looks like Figure~\ref{fig:gaussim},
where each $\alpha_j$ is the image of a planar geodesic in the
boundary of $M_0$.  
We identify $M_0$ with the image of $M_0$ by $G_0$.
 
\begin{remark}
  If the ``catenoid cousin fence" exists, 
  its hyperbolic Gauss map is not singly periodic.
  Indeed,
  if $G$ is singly periodic, 
  any fundamental pieces coincide mutually completely
  (not only are congruent).
  This contradiction implies that $G$ cannot be singly periodic.
\end{remark}

\begin{figure}
\begin{center}
\unitlength=1pt
\begin{picture}(403.00,122.00)(98.00,688.00)
\put(163.00,771.00){\line(1,0){10.00}}
\put(168.00,725.00){\line(1,0){11.00}}
\put(158.00,689.00){\line(1,0){8.00}}
\put(340.00,730.00){\line(1,0){160.00}}
\put(420.00,810.00){\line(0,-1){117.00}}
\put(419.00,790.00){\line(0,-1){39.00}}
\put(421.00,751.00){\line(0,1){30.00}}
\put(220.00,760.00){\vector(1,0){101.00}}
\put(160.00,725.00){\circle*{4}}
\put(154.00,740.00){\makebox(0,0)[cc]{$\alpha_1$}}
\put(132.00,760.00){\makebox(0,0)[cc]{$\alpha_2$}}
\put(129.00,732.00){\makebox(0,0)[cc]{$\alpha_3$}}
\put(121.00,720.00){\makebox(0,0)[cc]{$\alpha_3$}}
\put(126.00,699.00){\makebox(0,0)[cc]{$\alpha_4$}}
\put(178.00,708.00){\makebox(0,0)[cc]{$\frac{\pi}{2}$}}
\put(180.00,751.00){\makebox(0,0)[cc]{$\frac{\pi}{2}+\delta$}}
\put(165.00,730.00){\makebox(0,0)[cc]{$p$}}
\put(192.00,787.00){\makebox(0,0)[cc]{$S^2$}}
\put(266.00,747.00){\makebox(0,0)[cc]{stereographic}}
\put(266.00,737.00){\makebox(0,0)[cc]{projection from}}
\put(266.00,727.00){\makebox(0,0)[cc]{the point antipodal}}
\put(266.00,717.00){\makebox(0,0)[cc]{to $p$}}
\put(424.00,725.00){\makebox(0,0)[cc]{$p$}}
\put(397.00,725.00){\makebox(0,0)[cc]{$\alpha_1$}}
\put(468.00,762.00){\makebox(0,0)[cc]{$\alpha_2$}}
\put(427.00,760.00){\makebox(0,0)[cc]{$\alpha_3$}}
\put(414.00,767.00){\makebox(0,0)[cc]{$\alpha_3$}}
\put(371.00,785.00){\makebox(0,0)[cc]{$\alpha_4$}}
\put(460.00,797.00){\makebox(0,0)[cc]{complex plane}}
\put(471.00,722.00){\makebox(0,0)[cc]{$\tan(\frac{\pi}{4}+\frac{\delta}{2})$}}
\put(360.00,723.00){\makebox(0,0)[cc]{-1}}
\bezier100(149.00,790.00)(98.00,789.00)(100.00,740.00)
\bezier101(149.00,790.00)(201.00,789.00)(200.00,740.00)
\bezier101(200.00,740.00)(201.00,690.00)(150.00,690.00)
\bezier100(150.00,690.00)(101.00,689.00)(100.00,740.00)
\bezier57(149.00,790.00)(162.00,784.00)(160.00,740.00)
\bezier55(160.00,740.00)(161.00,696.00)(150.00,690.00)
\bezier87(112.00,734.00)(113.00,775.00)(159.00,771.00)
\bezier36(112.00,732.00)(106.00,730.00)(136.00,726.00)
\bezier47(100.00,740.00)(104.00,726.00)(137.00,725.00)
\bezier38(164.00,688.00)(175.00,710.00)(174.00,724.00)
\bezier48(175.00,725.00)(174.00,761.00)(167.00,771.00)
\bezier121(419.00,790.00)(359.00,791.00)(360.00,730.00)
\bezier96(422.00,781.00)(468.00,780.00)(470.00,730.00)
\end{picture}
\end{center}
        \caption{Modified Fundamental Region $M_{\delta,0}$.}%
        \label{fig:fenceregion}
\end{figure}
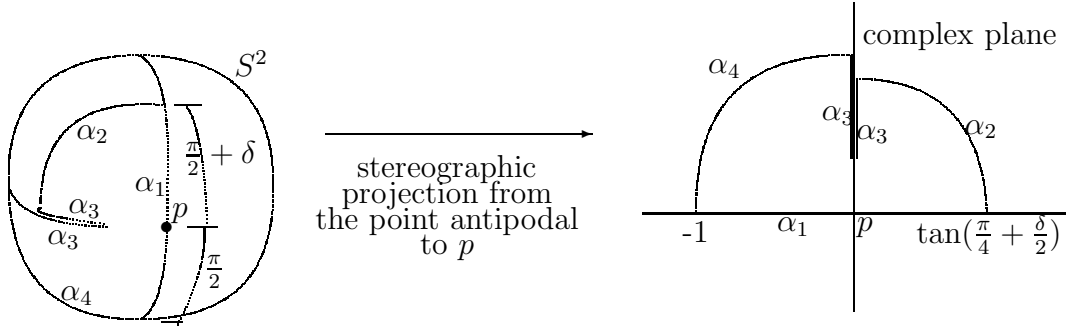

For each number $\delta\in(-\pi/2,\pi/2)$,  
we define a region $M_{\delta,0}$ on $S^2$ as in 
Figure~\ref{fig:fenceregion},
where $\alpha_1$ and $\alpha_3$ are segments 
of a great circles through $P$,
$\alpha_4$ is a segment of a great circle centered at $P$,
and $\alpha_2$ is one-fourth of a {\it smaller\/} circle centered
at $P$ with radius $\pi/2+\delta$.

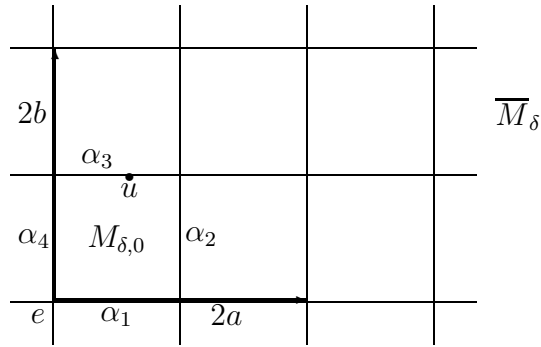
\begin{figure}
\begin{center}
\unitlength=0.8pt
\begin{picture}(240.00,160.00)(100.00,630.00)
\put(120.00,790.00){\line(0,-1){160.00}}
\put(180.00,790.00){\line(0,-1){160.00}}
\put(240.00,790.00){\line(0,-1){160.00}}
\put(300.00,790.00){\line(0,-1){160.00}}
\put(100.00,770.00){\line(1,0){220.00}}
\put(100.00,710.00){\line(1,0){220.00}}
\put(100.00,650.00){\line(1,0){220.00}}
\put(120.00,651.00){\vector(1,0){120.00}}
\put(121.00,650.00){\vector(0,1){120.00}}
\put(156.00,709.00){\circle*{4}}
\put(156.00,703.00){\makebox(0,0)[cc]{$u$}}
\put(113.00,643.00){\makebox(0,0)[cc]{$e$}}
\put(150.00,643.00){\makebox(0,0)[cc]{$\alpha_1$}}
\put(190.00,681.00){\makebox(0,0)[cc]{$\alpha_2$}}
\put(141.00,718.00){\makebox(0,0)[cc]{$\alpha_3$}}
\put(111.00,680.00){\makebox(0,0)[cc]{$\alpha_4$}}
\put(202.00,644.00){\makebox(0,0)[cc]{$2a$}}
\put(110.00,740.00){\makebox(0,0)[cc]{$2b$}}
\put(340.00,739.00){\makebox(0,0)[cc]{$\overline{M}_\delta$}}
\put(150.00,679.00){\makebox(0,0)[cc]{$M_{\delta,0}$}}
\end{picture}
\end{center}
        \caption{The Riemann surface $\widetilde M_{\delta}$.}%
        \label{fig:fenceriem}
\end{figure}

The boundary curves of $M_{\delta,0}$ are circles,
so we reflect $M_{\delta,0}$ 
infinitely often by inversions with respect  
to the boundaries of it and of its copies.
Then we get an abstract Riemann surface ${\overline M}_{\delta}$
which is biholomorphic to $\C$.
We identify ${\overline M}_{\delta}$ with $\C$ and 
take $\vect{a}$ and $\vect{b}$ in $\C$ as in Figure~\ref{fig:fenceriem}.
We denote
\begin{equation}\label{eq:fenceendumb}    
  \begin{aligned}
       E&=e+2\Z\vect{a}+2\Z\vect{b},\\
       U&=(u+2\Z\vect{a}+2\Z\vect{b})\cup (-u+2\Z\vect{a}+2\Z\vect{b}),
  \end{aligned}
\end{equation}         
and define
$M_{\delta}=\C\setminus E$.
We need to find a $2\vect{b}$-periodic 
\cmcc\ immersion of a cylinder 
$f\colon{} \overline{M}_{\delta}/(2\Z\vect{a})\to \hypc$        
which has ends at $E$ and umbilic points at $U$.
    
Let $G_{\delta,0}\colon{}M_{\delta,0}\to M_{\delta,0}\subset S^2$
be the {\it identity map}.
Since the boundary curves of $M_{\delta,0}$ are circles,
$G_{\delta,0}$ can be extended to the holomorphic map
$G_{\delta}\colon{}\C={\overline M}_{\delta}\to S^2=\riem$.

Let $\mu_j$ be a reflection on 
$\overline{M}_{\delta}$, which is the lift
of the reflection about $\alpha_j$.
By the definition of $G_{\delta}$, it satisfies
$\overline{G_{\delta}\circ\tilde\mu_j}=\sigma(\mu_j)^{-1}
   \cdot G_{\delta}$,
where 
\begin{equation}\label{eq:fencGsym2}
  \begin{alignedat}{2}
    \sigma(\mu_1)&=\id, \qquad
    &\sigma(\mu_2)&=
       \begin{pmatrix} 0   & 
         i\tan\left(({\pi}/{4})+({\delta}/{2})\right)\\
         i\cot\left(({\pi}/{4})+({\delta}/{2})\right)&0  
       \end{pmatrix}, \\ 
    \sigma(\mu_3)&= 
       \begin{pmatrix}
          i  & 0 \\  0 & -i 
       \end{pmatrix}, \qquad
    &\sigma(\mu_4)&=
       \begin{pmatrix} 0  & i \\  i & 0 \end{pmatrix}
   \end{alignedat}
\end{equation}            
(see Figure~\ref{fig:fenceriem}).   
We give the following holomorphic 2-differential
$Q_{\delta}$ as a  Hopf differential on $M_{\delta}$.
\begin{lemma}\label{lem:fenceQ}
  There exists a meromorphic $2$-differential 
  $Q_{\delta}$ which has the following properties.
  \begin{enumerate}
   \item   $Q_{\delta}$ is doubly periodic with
           respect to the periods generated by $2\vect{a}$ 
           and $2\vect{b}$.
   \item   $Q_{\delta}$ has poles of order $2$ at $E$,
           and is regular on $\C\setminus E$.
   \item   $Q_{\delta}$ has zeroes of order $1$ at $U$,
           and is non-zero on $\C\setminus U$.
   \item   $\overline{Q_{\delta}\circ\tilde\mu_j}=Q_{\delta}$.
  \end{enumerate}
  Moreover, $Q_{\delta}$ is unique up to a real constant 
  factor.  
\end{lemma}

\begin{pf}
  If we set $Q_{\delta}=q_\delta(z)dz^2$,
  then the proof reduces to finding a suitable elliptic 
  function $q_\delta(z)$, which can be done 
  by using elementary elliptic function theory. 
\end{pf}

Hence, we can choose $G_{\delta}$ and $Q_{\delta}$ continuously on $\delta$
so that 
  $ds_{\delta}^2=
      (1+|G_{\delta}|^2)^2 \omega_{\delta}\bar\omega_{\delta}$
is positive definite,     
where $\omega_{\delta}=Q_{\delta}/G_{\delta}$.

\begin{theorem}\label{thm:fence}
  For a sufficiently small $c$, there exist a number $\delta$ and 
  a $2\vect{b}$-periodic complete conformal \cmcc\ immersion
  $f\colon{}\overline{M}_{\delta}/(2\Z\vect{a})\to \hypc$
  which has hyperbolic Gauss map $G_\delta$ and Hopf
  differential $Q_\delta$.
\end{theorem}                    

By the same argument as in the proof of Proposition~\ref{prop:JMStep2},
we have the following:
\begin{lemma}\label{lem:fenceStep1}
  There exists a two-parameter family 
  $\{F_{c,\delta}\}_{|c|<\varepsilon,|\delta|<\delta_0}$
  of null holomorphic immersions into $SL(2,\C)$        
  with the following properties:
  \begin{enumerate}
    \item \label{item:fence1}
          $f_{c,\delta}=F_{c,\delta}{F_{c,\delta}^*}\in
           \D_{\overline{M}_{\delta}/2\Z\vect{a}}^{(c)}
           (G_{\delta},Q_{\delta})$ 
                for each $(c,\delta)$, where $c\neq 0$\rom.
    \item \label{item:fence2}
          $F_{c,\delta}$ is smooth in $c$ and continuous in $\delta$.
    \item \label{item:fence3}
          $\lim_{c\to 0}F_{c,\delta}=\id$.
    \item \label{item:fence4}
          $\hat\rho_{F_{c,\delta}}(\tilde\mu_1)=\id$,  and
          \[
              \hat\rho_{F_{c,\delta}}(\tilde\mu_3)=
                \begin{pmatrix}
                  \xi & 0  \\
                  0 & \xi^{-1}
              \end{pmatrix},\qquad
              \hat\rho_{F_{c,\delta}}(\tilde\mu_4)=
                \begin{pmatrix}
                  p &   i\beta\\
                  i\beta  &   \bar p
                \end{pmatrix},\qquad                       
              \hat\rho_{F_{c,\delta}}(\tilde\mu_2)=
                \begin{pmatrix}
                  q &   i\gamma_1\\
                  i\gamma_2  &   -\bar q
                \end{pmatrix},
           \]                      
  \end{enumerate}
  where $\xi=\xi(c,\delta)$, $p=p(c,\delta)$
  and $q=q(c,\delta)$  are  complex-valued functions
  such that $|\xi|=1$,
  and $\beta$, $\gamma_1$ and $\gamma_2$ are real-valued functions 
  of $c$ and $\delta$.
  They satisfy
  \begin{gather*}
    p(0,\delta)=0, \qquad  q(0,\delta)=0, \qquad
    \beta(0,\delta)=1, \qquad  \xi(0,\delta)=i, \\
    \gamma_1(0,\delta)=
         \tan\left(\frac{\pi}{4}+\frac{\delta}{2}\right), \qquad
    \gamma_2(0,\delta)=
         \cot\left(\frac{\pi}{4}+\frac{\delta}{2}\right).
  \end{gather*}                            
\end{lemma}
\begin{pf*}{\it Proof of Theorem~\ref{thm:fence}}
  In Lemma~\ref{lem:fenceStep1}, let 
  $\phi(c,\delta)=\gamma_1(c,\delta)/\gamma_2(c,\delta)$.
  So $\phi(0,\delta)>1$ if $\delta>0$,
  and $\phi(0,\delta)<1$ if $\delta<0$.
  Then, by the continuity of $\gamma_j$, there exist $\delta_+$ and $\delta_-$
  such that $\phi(c,\delta_+)>1$ and $\phi(c,\delta_-)<1$ hold for a 
  sufficiently small $c$. Hence there exists $\delta$ such that
  $\phi(c,\delta)=1$.  For such $c$ and $\delta$,
  all matrices in Lemma~\ref{lem:fenceStep1} are in $SU(2)$.
\end{pf*}
\begin{remark}  
  The same argument can be applied to the
  Jorge-Meeks $n$-oid fence  (see Section~4 in \cite{Rossman}), and this
produces
  a Jorge-Meeks $n$-oid  fence cousin. 
\end{remark}                

\section{Appendix}\label{appendix}

Let $\Gamma$ be a subgroup of $SU(2)$.
In this appendix, we prove a property of
a set of groups conjugate to $\Gamma$ in $SL(2,\C)$ defined by
\[
  C_{\Gamma}:=
  \{  \sigma\in SL(2,{\bf C})\,;\,
      \sigma\Gamma \sigma^{-1}\subset SU(2)
  \}.
\]
The authors  wish to thank Hiroyuki Tasaki for 
valuable comments on the first draft of the appendix.

If $\sigma\in C_{\Gamma}$, it is obvious that
$a\sigma\in C_{\Gamma}$
for all $a\in SU(2).$
So if we consider the quotient space
\[
   I_{\Gamma}:=C_{\Gamma}/SU(2),
\]
the structure of the set $C_{\Gamma}$ is completely determined.
Define a map $\tilde\phi:C_{\Gamma}\to {\cal H}^3$ 
by
\[
   \tilde\phi(\sigma):=\sigma^*\,\sigma,
\]
where ${\cal H}^3$ is the hyperbolic 3-space defined by
${\cal H}^3:=\{aa^*\,;\, a\in SL(2,\C)\}$.
Then it induces an
injective map $\phi: I_{\Gamma}\to {\cal H}^3$
such that
$\phi\circ \pi=\tilde \phi$, where $\pi:C_{\Gamma}\to I_{\Gamma}$
is the canonical projection.
So we can identify $I_{\Gamma}$ with a subset
$\tilde\phi(I_{\Gamma})=\tilde\phi(C_{\Gamma})$
of the hyperbolic 3-space ${\cal H}^3$.
The following assertion holds.

\begin{applemma}
  The subset $\tilde\phi(I_{\Gamma})$ is a point, a geodesic line,
  or all ${\cal H}^3$.
\end{applemma}
\begin{appremark}
Theorem \ref{thm:red} in Section 3 
is directly obtained if we set $\Gamma:=\rho(\pi_1(M))$,
where $\rho$ is the representation defined in Section 3.
Indeed, the set $I_{M}^{(c)}(G,Q)$ defined in Section 3
coincides with the set $\tilde\phi(I_{\rho(\pi_1(M))})$.  
\end{appremark}

\begin{pf}
  For each $\gamma\in \Gamma$, we set
  \[
      C_\gamma:=\{c\in SL(2,{\bf C})\,;\,
                  \sigma \gamma \sigma^{-1}\in SU(2)
        \}.
  \]
  Then we have
  \begin{equation}\label{app:1}
      C_\Gamma:=\bigcap_{\gamma\in \Gamma} C_\gamma.
  \end{equation}
  The condition $\sigma\gamma \sigma^{-1}\in SU(2)$
  is rewritten as $\sigma^*\,\sigma\,\gamma=\gamma \sigma^*\,\sigma$.
  So we have
  \begin{equation}\label{app:2}
     \tilde \phi (C_\gamma)=\H^3 \cap Z_\gamma,
  \end{equation}
  where $Z_\gamma$ is the center of $\gamma\in \Gamma$.

  Assume $\gamma\neq \pm \id$.
  If $\gamma$ is a diagonal matrix, it can easily be checked that
  $Z_{\gamma}$ consists of diagonal matrices in $SL(2,\C)$.
  Since any $\gamma\in\Gamma$ can be diagonalized by a matrix in 
  $SU(2)$, we have $Z_{\gamma}=\{\exp( zT)\,;\,z\in\C\}$, where
  $T\in\Su(2)$ is chosen so that $\gamma=\exp(T)$.
  ($\Su(2)$ is the Lie algebra of $SU(2)$.)  
  Hence we have 
  \begin{equation}\label{app:3}
     \tilde\phi(C_{\gamma})=\H^3\cap Z_{\gamma}=
           \exp\left(i\R T\right),
  \end{equation}
  because $\exp\left(i\Su(2)\right)=\H^3$.

  Suppose now that $\Gamma$ is not abelian.
  Then there exist $\gamma$, $\gamma'\in \Gamma$
  such that $\gamma\gamma'\ne \gamma'\gamma$.
  Set $\gamma=\exp(T)$ and $\gamma'=\exp(T')$,
  where $T,T'\in \Su(2)$.
  Then we have $i {\R}T\cap i {\R}T'=\{0\}$.
  It is well-known that the restriction of the exponential map
  $\exp|_{i\Su(2)}\colon{}i\Su(2)\to\H^3$ is bijective.
  Hence we have
  \[
      \tilde\phi(C_\gamma)\cap\tilde\phi( C_{\gamma'})=
      \exp\left(i{\R}T\right)\cap 
      \exp\left(i{\R}T'\right)=\{\id\}.
  \]
  By (\ref{app:1}), (\ref{app:2}) and (\ref{app:3}),
  we have
  \[
      \tilde\phi(I_{\Gamma})=\{\id\} \qquad 
      \mbox{(if $\Gamma$ is not abelian)}.
  \]
  Next we consider the case $\Gamma$ is abelian.
  If $\Gamma\subset \{\pm\id\}$, then obviously
  \[
      \tilde\phi(I_{\Gamma})=\H^3 \qquad 
               \mbox{(if $\Gamma\subset \{\pm \id\}$)}.
  \]
  Suppose $\Gamma\not\subset \{\pm\id\}$. Then 
  there exists $\gamma\in \Gamma$ such that $\gamma\ne \pm \id$.
  We set $\gamma=\exp{T}$\,\,($T\in \Su(2)$).
  Since $\exp(\R T)$ is a maximal abelian subgroup
  containing $\gamma$, we have $\Gamma \subset \exp(\R T)$.
  Then by (\ref{app:3}), we have
  \[
      \tilde\phi(I_{\Gamma})=\exp(i\R T)\qquad
      \mbox{(if $\Gamma\not\subset \{\pm \id\}$ is abelian)}.
      \qquad
      \Box
  \]
\renewcommand{\qed}{}
\end{pf}

\setlength{\parindent}{0in}
\small
\vspace{0.5in}
\sc 
Wayne Rossman \par
Graduate School of Mathematics \par
Kyushu University \par
Fukuoka 812-81\par
Japan.\par
\it E-mail \normalshape: 
wayne@@math.kyushu-u.ac.jp

\vspace{0.2in}

\sc
Masaaki Umehara\par
Department of Mathematics \par
Graduate School of Science\par
Osaka University \par
Toyonaka, Osaka 560 \par
Japan.\par
\it E-mail \normalshape:
umehara@@math.wani.osaka-u.ac.jp

\vspace{0.2in}

\sc
Kotaro Yamada \par
Department of Mathematics \par
Faculty of Science \par
Kumamoto University \par
Kumamoto 860\par
Japan.\par
\it E-mail \normalshape: 
kotaro@@gpo.kumamoto-u.ac.jp

\end{document}